\documentclass{amsart}
\usepackage{amssymb, url, algorithm, algorithmic, multirow}

\newcommand{\ai}{\mathfrak{a}}
\newcommand{\ah}{\hat{\alpha}}
\newcommand{\bi}{\mathfrak{b}}
\newcommand{\C}{\mathbb{C}}
\newcommand{\ch}{\operatorname{char}}
\newcommand{\D}{\mathcal{D}}
\newcommand{\dv}{\operatorname{div}}
\newcommand{\f}{\mathfrak{f}}
\newcommand{\g}{\mathfrak{g}}
\newcommand{\F}{\mathbb{F}}
\newcommand{\Fq}{\F_q}
\newcommand{\HH}{\mathcal{H}}
\newcommand{\I}{\mathcal{I}}
\newcommand{\ideal}[1]{\left\langle {#1} \right\rangle }
\newcommand{\J}{\mathcal{J}}
\newcommand{\kk}{\mathfrak{k}}

\newcommand{\N}{\mathbb{N}}
\newcommand{\OO}{\mathcal{O}}
\newcommand{\p}{\mathfrak{p}}
\newcommand{\PP}{\mathcal{P}}

\newcommand{\Pf}{\mathfrak{P}}
\newcommand{\ph}{\hat{\psi}}
\newcommand{\ppmod}[1]{\mbox{\ $\left(\operatorname{mod}\ {#1}\right)$}}
\newcommand{\q}{\mathfrak{q}}
\newcommand{\R}{\mathcal{R}}

\newcommand{\sig}{\operatorname{sig}}
\newcommand{\sgn}{\operatorname{sgn}}
\newcommand{\st}{\mid}
\newcommand{\supp}{\operatorname{supp}}
\newcommand{\SC}{\mathcal{S}}
\newcommand{\tk}{\mathfrak{t}}
\newcommand{\wk}{\mathfrak{w}}
\newcommand{\Z}{\mathbb{Z}}

\newtheorem{theorem}{Theorem}[section]

\newtheorem{proposition}[theorem]{Proposition}

\newtheorem{heuristic}[theorem]{Heuristic}

\theoremstyle{definition}

\theoremstyle{remark}

\numberwithin{algorithm}{section}
\numberwithin{equation}{section}

\begin{document}

\title{Class number and regulator computation in cubic function fields}


%    author one information
%\author[E. Landquist]{Eric Landquist}
%\address{Carl von Ossietzky Universit{\"{a}}t Oldenburg,
%Institut f{\"{u}}r Mathematik,
%Ammerl{\"{a}}nder Heerstra{\ss}e 114-118,
%26129 Oldenburg, Germany}
%\email{eric.landquist@uni-oldenburg.de}
%\curraddr{}
%\thanks{}

\author[E. Landquist]{Eric Landquist}
\address{Kutztown University,
Department of Mathematics, Kutztown, PA 19530} \email{elandqui@kutztown.edu}
%\thanks{}

%    author two information
\author[R. Scheidler]{Renate Scheidler}
\address{University of Calgary,
Department of Mathematics and Statistics, 2500 University Drive NW, Calgary, Alberta T2N 1N4, Canada} \email{rscheidl@ucalgary.ca}
\thanks{Research supported by NSERC of Canada.}

%    author three information
\author[A. Stein]{Andreas Stein}
\address{Carl von Ossietzky Universit{\"{a}}t Oldenburg,
Institut f{\"{u}}r Mathematik, D-26111 Oldenburg, Germany} \email{andreas.stein1@uni-oldenburg.de}
%\thanks{}

%    \subjclass is required.
\subjclass[2000]{11R16, 11R58, 11R29, 11Y16, 11Y40, 14H05, 14Q05} % 11R65, 11M38, 11G20

\date{}

\dedicatory{}

%    Abstract is required.
\begin{abstract}
We present computational results
%on a method of Scheidler and Stein
on
%computing
the divisor class number and the regulator of a cubic function field over a large base field. The underlying method is based on approximations
of the Euler product representation of the zeta function of such a field. We give details on the implementation for purely cubic function
fields of signatures $(3,1)$ and $(1, 1; 1, 2)$, operating in the ideal class group and infrastructure of the function field, respectively. Our
implementation provides numerical evidence of the computational effectiveness of this algorithm. With the exception of special cases, such as
purely cubic function fields defined by superelliptic curves, the examples provided are the largest divisor class numbers and regulators ever computed
for a cubic function field over a large prime field. The ideas underlying the optimization of the class number algorithm can in turn be used to
analyze the distribution of the zeros of the function field's zeta function. We provide a variety of data on a certain distribution of the
divisor class number that verify heuristics by Katz and Sarnak on the distribution of the zeroes of the zeta function.
\end{abstract}

\maketitle


\section{Introduction and Motivation}
\label{S:intro} One of the more difficult problems in arithmetic geometry is the computation of the order of the group of rational points on the
Jacobian of an algebraic curve over a finite field. In this paper, we give results on the application and optimization of a method of Scheidler
and Stein \cite{ss07, ss08} to compute the order of the Jacobian of a purely cubic function field over a large base field. We describe details
of this implementation for the cases in which the function field has signature $(3,1)$ and $(1,1;1,2)$, operating in the ideal class group and
 infrastructure of the function field, respectively. We also provide numerical results for both scenarios.

%For a complete background on this topic, we refer to \cite{ss07,ss08}, where the interested reader will find an extensive overview of related work.

 In general, determining the divisor class number $h$ of an algebraic function field $K/\Fq(x)$ 
of genus $g$ with $q=p^n$ large is considered to be a computationally difficult problem. 
There are several methods available to compute this class number. Some of these methods are 
general, while others apply only to specific curves. Here, we only highlight the literature 
that is most closely related to our context. 
Kedlaya's $p$-adic algorithm, using Monsky-Washnitzer cohomology, computes the zeta function of 
a hyperelliptic curve over a finite field of odd characteristic \cite{k, ke, kedlaya2004} and
requires $\tilde{O}(pn^3g^4)$ bit operations, for fixed $g$ and $n$. 
This method has been generalized 
to superelliptic curves \cite{gg, lauder}, 
$C_{ab}$ curves \cite{dv06}, and nondegenerate curves \cite{cdv}. 
%Harvey adapted Kedlaya's algorithm and has a time complexity of
%$\tilde{O}\left(p^{0.5}n^{3.5}g^{5.5}+\log(p)n^5g^8\right)$ bit operations, for fixed $g$ and $n$, 
%to compute the zeta function of hyperelliptic curves, including a $46$-digit class number of
%a genus $3$ hyperelliptic curve in 40 hours and a $53$-digit class number of a genus $4$
%curve in 45 hours.
% g = 3, 46-digit (~2^{150}) class number in 40 hours.
% g = 4, 53-digit (~2^{176}) class number in 45 hours. 
Minzlaff adapted Kedlaya's algorithm to compute the zeta function of a superelliptic curve with
$\tilde{O}\left(p^{0.5}n^{3.5}g^{5.5}+\log(p)n^5g^8\right)$ bit operations, for fixed $g$ and $n$,   
\cite{minzlaff} and found
the zeta function of a Picard curve over $q = 2^{36}+31$ in 9.4 hours.
While not computed explicitly, the class number of this curve would have been $33$ digits.
% Used a Picard curve, computed zeta function with $q = 2^{36}+31$ in 9.4 hours.
% 33-digit, (~2^{108}) class number in 9.4 hours.
For larger genera, index calculus methods have been developed to compute discrete logarithms 
(and hence group orders) in class groups arising from plane curves in expected time 
$\tilde{O}\left( q^{2-2/g} \right)$ bit operations, for fixed genus, as $q \to \infty$ 
\cite{gttd}, and 
from plane curves of small degree in expected time $\tilde{O}\left( q^{2-2/(g-1)} \right)$, 
for fixed genus, as $q \to \infty$ \cite{diem06}. In \cite{dt08}, the latter algorithm was 
tested on a Koblitz $C_{3,4}$ curve (of genus $3$) with $q = 2^{31}$ and completed in 
9.3 CPU days. Each of these methods requires
the place at infinity to ramify completely.

The algorithm of Hess \cite{hess-phd}, on the other hand, applies to any global 
function field, provided there is at least one infinite place of degree $1$ 
and computes the structure of the divisor class group in expected time 
$O\left(\exp\left(\sqrt{2(2g\log(q))\log(2g\log(q))}\right) \right)$ bit operations, 
as $g \to \infty$. This algorithm is implemented in MAGMA and later we will compare 
the efficiency of this method with ours.

A straightforward
method is to search for $h$ in the Hasse-Weil interval $[(\sqrt{q}-1)^{2g},\,(\sqrt{q}+1)^{2g}]$ using Shanks' Baby Step-Giant Step
\cite{shanksbsgs} or Pollard's Kangaroo \cite{pollardkangaroo} method. In this way, $h$ can  be computed deterministically or heuristically,
respectively, using $O\left(q^{(2g-1)/4}\right)$ group operations as $q \to \infty$.

Stein and Williams \cite{sw} applied techniques used by
Lenstra \cite{len} and Schoof \cite{schoof82} in quadratic number fields to real hyperelliptic function fields to narrow the search
space. The basic idea is to find an approximation $E$ of $h$ and value $U$ such that $|h-E|<U$. The new interval $[E-U,\,E+U]$ is much smaller
than the Hasse-Weil interval and in practice much better than expected. As a result, the Stein-Williams algorithm finds the divisor class number
and regulator using $O\left(q^{[(2g-1)/5] + \varepsilon(g)}\right)$ infrastructure operations, where $-1/4\leq \varepsilon(g)
\leq 1/2$. The method was generalized to arbitrary hyperelliptic functon fields and improved by Stein and Teske \cite{st02,st02b,st05}, who also
applied Pollard's Kangaroo algorithm to this setting to compute the $29$-digit class number and regulator of a real hyperelliptic function field
of genus $3$.

The algorithm of \cite{sw} was generalized to cubic function fields in \cite{ss07} and to arbitrary function fields in \cite{ss08}. In this paper,
we provide an implementation and mumerical examples for purely cubic function fields of signatures $(3,1)$ and $(1,1;1,2)$.
For this implementation, one requires
efficient arithmetic as well as effective criteria for determining how a place of $\Fq(x)$ splits in $K$. These and other algorithmic details can
be found in \cite{ericphd}. We show how to explicitly compute an estimate $E$ of $h$ and an upper bound $U$ on the error $|h-E|$ based on the
methods of \cite{ss07,ss08}. Furthermore, we provide experimental results on the distribution of $h$ in the interval $[E-U,\,E+U]$ and show how
to optimize Pollard's Kangaroo algorithm on this interval, making improvements to its application to the infrastructure of a purely cubic function
field of signature $(1,1;1,2)$. The improved method is applied to compute divisor class numbers up to $28$ and $25$ digits of function fields
of genus $3$ and $4$, respectively. For the signature $(1,1;1,2)$ examples, we extracted the regulators as well. These are the largest known
class numbers and regulators ever computed for a cubic function field over a large base field, with the exception of those computed by Bauer,
Teske, and Weng in \cite{btw}.

A brief discussion of the method of \cite{btw}, and its more memory-efficient variant of \cite{weng}, is appropriate here. This technique is
specific to function fields defined by a Picard curve, for which it determines the divisor class number using $O\left(\sqrt{q}\right)$ Jacobian
operations. It has generated a divisor class number as large as $55$ digits, the largest known class number of any cubic function field.
However, this technique is restricted to Picard curves, whereas ours is a general purpose algorithm, and our implementation applies to any purely
cubic function field of signature $(3,1)$ or $(1,1;1,2)$. Picard fields are a special case of our setting, namely purely cubic function fields
of signature $(3,1)$ for which the curve under consideration must be nonsingular and of genus 3.

The remainder of this paper is organized as follows. We give an overview of cubic function fields, including their class group and (in the
signature $(1,1;1,2)$ case) their infrastructure, in Section \ref{S:background}. In Section \ref{S:kangaroo}, we describe Pollard's Kangaroo
method, outlining its specific application to the ideal class group of a purely cubic function field of signature $(3,1)$. We then adapt the
Kangaroo method to the infrastructure of a purely cubic function field of signature $(1,1;1,2)$ in Section \ref{S:rooinf}. In Section
\ref{S:EU}, we review the main results of \cite{ss07} to compute the divisor class number of a cubic function field and explain the details
of our implementation in Section \ref{S:details}. Finally, we provide numerical results obtained by our implementation in Section \ref{S:results}
and conclude with open problems and areas for further research in Section \ref{S:conclusions}.

\section{Cubic Function Fields --- Class Group and Infrastructure}
\label{S:background}

For a general introduction to function fields, we direct the reader to \cite{h, sti, rosen}. Details on cubic function fields and their arithmetic
can be found in \cite{ss00,s01,s04,b,ss07,ericphd}. Here, we merely highlight the material that is required for our context. Let $\Fq$ be a finite
field and $\Fq(x)$ the field of rational functions in $x$ over $\Fq$. Throughout this paper, we assume that $\ch(\Fq) \geq 5$. A {\em cubic function
field} is a separable extension $K/\Fq(x)$ of degree $3$ with full constant field $\Fq$; we write $K = \Fq(C)$ with $C: f(x,\,Y)=0$, where
$f(x,\,Y)\in\Fq[x,\,Y]$ is an absolutely irreducible monic polynomial of degree $3$ in $Y$. Note that we do not require the curve $f(x,\,Y) = 0$
to be nonsingular. If $C: Y^3 = F$ with $F \in \Fq[x]$ cube-free, then $K$ is called a {\em purely} cubic function field.  In this case, if we
write $F = GH^2$ with $G,\,H \in \Fq[x]$, $\gcd(G,\, H) = 1$, and both $G$ and $H$ are square-free, then the genus $g$ of $K$ is given by $g =
\deg(GH)-2$ if  $3 \mid\deg(F)$ and $g = \deg(GH)-1$ if $3 \nmid\deg(F)$.

\subsection{Divisors and Ideals}
\label{S:divisors}

Let $\D$ denote the group of divisors of $K$ defined over $\Fq$, $\D_0$ the subgroup of divisors of degree $0$ defined over $\Fq$, and $\PP$ the
subgroup of principal divisors defined over $\Fq$. Then the {\em $($degree~$0)$ divisor class group} (or {\em Jacobian}) of $K$ is the quotient
group $\J = \D_0/\PP$ and its order $h = |\J|$ is the {\em $($degree~$0)$ divisor class number} of $K$. Let $S$ be the collection of places of
$K$ lying above the place at infinity of $\Fq(x)$, $\D_0^S = \{D\in\D_0\st \supp(D) \subseteq S\}$, and $\PP^S$ the subgroup of principal divisors
in $\D_0^S$. Then the order $R_x$ of the quotient group $\D_0^S/\PP^S$ is the {\em regulator} of $K$. Finally, let $\D_S = \{D\in\D\st
\supp(D)\cap S = \emptyset\}$ and $\PP_S = \PP\cap\D_S$. Then every divisor $D\in\D$ can be uniquely written in the form $D = D_S + D^S$ with
$D_S\in\D_S$ and $D^S\in\D^S$.

The {\em maximal order} of $K/\Fq(x)$ (or {\em coordinate ring} of $C$) is the integral closure of $\Fq[x]$ in $K$ and is denoted $\OO_x$.
%By Dirichlet's Unit Theorem, the {\em unit rank} of $K$, i.e.\ the rank of the unit group $\OO_x^*$ of $\OO_x$ is $|S|-1$.
Let $\I = \I(\OO_x)$ be the group of fractional ideals of $\OO_x$ and $\HH =\HH(\OO_x)$ the subgroup of principal fractional ideals.
The {\em ideal class group} of $K$ is the quotient group $Cl(\OO_x) = \I/\HH$, and its order $h_x = |Cl(\OO_x)|$ is the {\em ideal
class number} of $K$. Let $f_x$ be the greatest common divisor of the inertia degrees of all the places in $S$. By Schmidt \cite{schmidt},
there is an exact sequence
$$ (0) \to \D_0^S/\PP^S \to \J \to Cl(\OO_x) \to \Z/f_x\Z \to (0) \enspace , $$
so that $f_xh = R_xh_x$.

There is a well-known isomorphism  $\Phi: \D_S \to \I$ given by $D \mapsto \{\alpha \in K^* \st \dv(\alpha)_S \geq D\}$ with inverse $\f \mapsto
\sum_{\p\notin S} m_{\p}\p$, where $\p$ denotes any finite place of $K$, $m_{\p} = \min\{v_{\p}(\alpha) \st \alpha \in \f \mbox{ non-zero}\}$, and
$v_{\p}$ is the normalized discrete valuation corresponding to $\p$. Moreover, $\Phi$ induces an isomorphism from $\D_S/\PP_S$ to $Cl(\OO_x)$.

If $S$ contains an infinite place $\infty_0$ of degree $1$, then $\Phi$ can be extended to an isomorphism
$$\Psi: \big \{ \, D \in \D_0 \st v_{\p}(D) = 0 \mbox{ for all } \p \in S \setminus \{ \infty_0\} \, \big \} \to \I(\OO_x)$$
via $\Psi(D_S - \deg(D_S)\infty_0) = \Phi(D_S)$, with inverse map $\Psi^{-1}(\f) = \Phi^{-1}(\f) - \deg(\f)\infty_0$.

\subsection{Splitting of Places}
\label{S:signatures}

Via $\Phi$, there is a one-to-one correspondence between the finite places of $K$ and the prime ideals of $\OO_x$. If $\Pf$ is a finite place
of $\Fq(x)$, then $\Phi(\Pf) = \ideal{P}$, the principal ideal generated by some irreducible polynomial $P\in\Fq[x]$, and the splitting behavior
of $\Pf$ in $K$ is identical to the splitting behavior of $\ideal{P}$ in $\OO_x$. Therefore, we may characterize the splitting behavior of the
finite places of a purely cubic function field by considering the following criteria.

\begin{theorem}[Theorem 3.1 of \cite{s01}] \label{thm:primes}
Let $K = \Fq(C)$ be a purely cubic function field with $C:Y^3 = F \in \Fq[x]$ cubefree and $\ch(K)\geq 5$.
If $P \in \Fq[x]$ is an irreducible polynomial, then the principal ideal $\ideal{P}$ splits into prime ideals in $\OO_x$ as follows:
\begin{enumerate}
\item If $P\mid F$, then $\ideal{P} = \p^3$.
\item If $P\nmid F$ and $q^{\deg(P)} \equiv 2 \ppmod{3}$, then $\ideal{P} = \p\q$.
\item If $P\nmid F$, $F$ is a cube modulo $P$, and $q^{\deg(P)} \equiv 1 \ppmod{3}$, then $\ideal{P} = \p\p'\p''$.
\item If $P\nmid F$ and $F$ is not a cube modulo $P$, then $\ideal{P} = \p$.
\end{enumerate}
\end{theorem}

An algorithmic realization of this theorem was given in \cite{ss07}.

We also require the splitting behavior of the infinite place $\infty$ of $\Fq(x)$ in $K$. Write $S = \{ \infty_0, \infty_1, \ldots , \infty_r \}$.
Let $e_i = e(\infty_i|\infty)$ be the ramification index and $f_i = f(\infty_i|\infty)$ the inertia degree of $\infty_i$ over $\infty$ for
$0 \leq i \leq r$. The {\em signature} of $K/\Fq(x)$ is the $(2r+2)$-tuple  $\sig(K) = (e_0, f_0; \ldots ; e_r, f_r)$  and is given as follows:

\begin{theorem}[Theorem 2.1 of \cite{ss00}] \label{cor:signatures}
Let $K = \Fq(C)$ be a purely cubic function field with $C:Y^3 = F \in \Fq[x]$ cube-free and $\ch(K)\geq 5$. If $\sgn(F)$ denotes the leading
coefficient of $F$, then the signature of $K/\Fq(x)$ is given as follows:
\begin{enumerate}
\item If $3\nmid\deg(F)$, then $\sig(K) = (3,1)$.
\item If $3\mid\deg(F)$ and $q \equiv 2 \ppmod{3}$, then $\sig(K) = (1,1;1,2)$.
\item If $3\mid\deg(F)$, $\sgn(F) \in \Fq^3$, and $q \equiv 1 \ppmod{3}$, then $\sig(K) = (1,1;1,1;1,1)$.
\item If $3\mid\deg(F)$ and $\sgn(F) \notin \Fq^3$, then $\sig(K) = (1,3)$.
\end{enumerate}
\end{theorem}

Thus, the signature of a purely cubic function field can be determined quickly. We also note that whenever $K$ has an infinite place of degree
$1$, there is a straightforward normalization that replaces $C:Y^3=F$ by an isomorphic curve $C':Y^3=F'$ where $F'$ is monic, without changing
the coordinate ring or the signature. Hence, in our implementation, we only consider monic polynomials $F$.

For the remainder of this paper, we assume that $S$ contains an infinite place $\infty_0$ of degree $1$, so that $f_x=1$ and $h = R_xh_x$.
Thus, we will no longer consider cubic function fields of signature $(1,3)$. In fact, our later focus will be on purely cubic function fields
of the two signatures $(3,1)$ and $(1,1;1,2)$. If $\sig(K) = (3,1)$, then $R_x = 1$ so that $\J \cong Cl(\OO_x)$ and $h=h_x$. In this case, we may
conduct our computations to find $h$ in either $\J$ or $Cl(\OO_x)$. If $\sig(K) = (1,1;1,2)$, then $h=R_xh_x$, and we operate in a certain set of
divisors called the {\em infrastructure} of $K$. To that end, and in order to operate effectively and explicitly in the infrastructure and the
class group of $K$, we require the notion of a {\em distinguished} divisor.

\subsection{Distinguished Divisors and Infrastructure}
\label{S:distinguished}

Let $K$ be a cubic function field with an infinite place $\infty_0$ of degree $1$ and maximal order $\OO_x$. A divisor $D$ of $K$ is said to be
{\em finitely effective} if $D_S\geq 0$; that is, $v_{\p}(D) \geq 0$ for all finite places $\p$ of $K$. Following \cite{b,gps,ericphd}, a finitely
effective divisor $D$ is said to be {\em distinguished} if
\begin{enumerate}
\item $D$ is of the form $D = D_S - \deg(D_S)\infty_0$, and
\item if $E$ is any finitely effective divisor that is linearly equivalent to $D$, with $\deg\left(E_S\right) \leq \deg\left(D_S\right)$ and
$E^S \geq D^S$, then $D=E$.
\end{enumerate}
An ideal $\ai$ of $\OO_x$ is said to be distinguished if $\Psi^{-1}(\ai)$ is a distinguished divisor. Note that every distinguished divisor $D$ is
uniquely determined by its finite part $D_S$ and hence corresponds to a unique distinguished ideal $\Psi(D)$. Thus, arithmetic of distinguished
divisors is reduced to ideal arithmetic. In the context of this paper, if suffices to know that such an arithmetic exists, and that it is
efficient.

By \cite[Corollaries 5.2 and 5.3]{b}, if $\sig(K) = (3,1)$, then every divisor class of $K$ contains a unique distinguished representative.
Analogously, every ideal class of $K$ contains a unique distinguished ideal.  The {\em composition} $\ai*\bi$ of $\ai$ and $\bi$ is the distinguished
ideal equivalent to the (generally non-distinguished) ideal $\ai\bi$. Similarly, the composition of two distinguished divisors $D_1, D_2$ is the
unique distinguished divisor that is linearly equivalent to the sum $D_1+D_2$. Arithmetic in $Cl(\OO_x)$ and $\J$ can thus be efficiently conducted
via these unique distinguished representatives when $\sig(K) = (3,1)$. For details, see \cite{b} and \cite{ericphd}.

Assume now that $\sig(K) = (1, 1; 1, 2)$. In this case, not every divisor class contains a distinguished divisor; however, if such as divisor
exists in any given class, then uniqueness still holds. In this scenario, we operate in the (principal) infrastructure of $K$. A general
treatment of infrastructures in function field extensions of arbitrary degree can be found in \cite{felix08,felix08phd}. The cubic scenario was
first presented in \cite{ss00,s01}, and we use the divisor-theoretic description of \cite{ericphd} here. We only discuss the basic
infrastructure operations and the notion of distance, again omitting all arithmetic details.

The (finite) set
$$\R = \left\{ D\in\D_0  \st D \mbox{ is distinguished and $\Psi(D)$ is a principal ideal}  \right\}$$
is the {\em (principal) infrastructure} of of $K$.

Since $K$ has exactly one infinite place $\infty_0$ of degree~1, there exists a unique embedding of $K$ into the field of Laurent series in
$x^{-1}$; that is, into the completion of $K$ at $\infty_0$. This embedding establishes the notion of degree and sign (i.e.\ leading
coefficient) of elements in $K$. Note that $K$ has a fundamental unit $\epsilon$ with $\deg(\epsilon) = 2R_x$ that is unique up to sign. If
$D\in\R$ and $\ai = \Psi(D)$, then there is a unique function $\alpha \in K^*$ such that $\ai = \ideal{\alpha}$ and $0\leq \deg(\alpha) <
\deg(\epsilon) = 2R_x$. The {\em distance} of $D$ is defined as $\delta(D) = \deg(\alpha)$. Therefore, we may order the divisors in $\R$
increasingly by distance:
$$\R = \{D_0 = 0,\,D_1,\,\ldots,\,D_{l-1} \}\enspace,$$
with $\delta(D_i) < \delta(D_{i+1})$ for $0\leq i < l-1$.\
There are two main operations on $\R$, the {\em baby step} and {\em giant step} operations. A baby step maps $D_i$ to $D_{i+1}$ for $0\leq i <
l-1$ and $D_{l-1}$ to $D_0$. We write $bs(D_i) = D_{i+1}$. The giant step operation is analogous to composition. As before, if
$D_1,\,D_2\in\R$, then $D_1+D_2$ is generally not distinguished. However, there is a uniquely defined and efficiently computable
function $\psi\in K^*$ such that
\begin{equation}
\label{eq:gs}
D_1\oplus D_2 = \Psi^{-1}\left(\ideal{\psi}\Psi(D_1+D_2)\right) \in \R
\end{equation}
and $-2g \leq \deg(\psi) \leq 0$. We call $\oplus$ the giant step operation. Under $\oplus$, $\R$ is an Abelian group-like
structure, failing only associativity. However,
$$ \delta(D_1\oplus D_2) = \delta(D_1) + \delta(D_2) + \deg(\psi) \enspace , $$
so that $D_1\oplus D_2$ is
``close to'' $D_1+D_2$ in terms of distance. A third operation that will be required is the computation of the divisor {\em below} any integer
$n\in\N$ with $0\leq n < 2R_x$. This is the unique divisor $D(n)\in\R$ such that $\delta(D(n)) \leq n < \delta(bs(D(n)))$. Details on how to compute a
baby step, giant step, and $D(n)$ in $\R$ are given in \cite{s01,ericphd}. We have the following heuristics on baby steps and giant steps, which
are based on extensive experimental results and plausible theoretical assumptions.

\begin{heuristic}
\label{the-heuristic} Let $K$ be a cubic function field of signature $(1,1;1,2)$ and genus $g$.
\begin{enumerate}
\item If $D\in\R$, then with probability $1-1/q$, we have $\deg(D_S) = g$.
\item If $D_1,\,D_2\in\R$ and $\Psi(D_1\oplus D_2) = \ideal{\psi}\Psi(D_1)\Psi(D_2)$, then with probability $1-O(1/q)$, we have
  $$\ph(g)= \left\{\begin{array}{ll}
  -\lfloor g/3\rfloor & \mbox{ if } g \not\equiv 1 \ppmod{3} \ , \\
  -(g+2)/3 & \mbox{ if } g \equiv 1 \ppmod{3} \ , \end{array}\right. $$
where we set $\ph(g) = \deg(\psi)$.
\end{enumerate}
\end{heuristic}
Based on this heuristic, we can conclude that a baby step in $\R$ has length $2$ with probability $1-O(1/q)$. Note that this corresponds exactly
to the observations made in real hyperelliptic function fields (see \cite{st02,st02b,st05}).

\section{The Kangaroo Method in $Cl(\OO_x)$}
\label{S:kangaroo}

If we are given integers $E,\,U\in\N$ such that $h \in [E-U,\,E+U]$, then the Baby Step-Giant Step and Kangaroo methods may be optimized to
compute $h$ with a deterministic and probabilistic running time of $O\left(\sqrt{U}\right)$ group operations, respectively. The Baby Step-Giant
Step method is generally faster than the Kangaroo method, but it requires the storage of $O\left(\sqrt{U}\right)$ group elements and cannot be
parallelized efficiently. For larger computations, the Kangaroo method is preferable, since variants of this algorithm require very little
storage and can be parallelized, hence our reason for describing and using this method. Specifically, we describe the parallelized Kangaroo
method of van Oorschot and Wiener \cite{vow,st05} and explain important improvements that apply in particular to the problem of computing the
divisor class number of a purely cubic function field of signature $(3,1)$.  We always assume that $2(E-U) > E+U$. In this way, if we determine
$h_0 \in [E-U,\,E+U]$ to be a multiple of $h$, then $h_0/2 < (E+U)/2 < E-U$. Thus, $h_0 = h$. It is important to note that this is not an
unreasonable restriction for the integers $E$ and $U$ produced by Scheidler and Stein's method in \cite{ss07}; only function fields over very
small base fields fail this criterion.

We now describe in detail a modification of the parallelized Kangaroo method using notation similar to that of \cite{st02,st05} for
hyperelliptic function fields. Let $m$ be the (even) number of available processors. The Kangaroo algorithm uses two {\em herds} of {\em
kangaroos}, a herd $\{T_1,\,\ldots,\,T_{m/2}\}$ of {\em tame} kangaroos, and a herd $\{W_1,\,\ldots,\,W_{m/2}\}$ of {\em wild} kangaroos. A
kangaroo is a sequence of distinguished ideals. The Kangaroo method requires a {\em collision} between a tame and a wild kangaroo to obtain $h$.
The tame kangaroos begin their jumps at distinct known points near $E$, and the wild kangaroos at points near $h$ whose location in the interval
is unknown, hence their respective names.

The idea of the algorithm is as follows. Let $\g$ be a distinguished ideal. Define a set of small (relative to $U$) random positive integers
$\{s_1,\,\ldots,\, s_{64}\}$, the {\em jump set} $J = \{\g^{s_1},\,\ldots,\,\g^{s_{64}}\}$, and a hash function $v:\I(\OO_x) \to \{1,\,\ldots,\,
64\}$. Initialize each tame kangaroo $T_i$ at a distinguished ideal $\tk_{0,i} \sim \g^{E + (i-1)\nu}$ for some small $\nu \in \Z$ and $1\leq
i\leq m/2$, and each wild kangaroo $W_j$ at a distinguished ideal $\wk_{0,j} \sim \g^{(j-1)\nu}$ for $1\leq j\leq m/2$. The kangaroos jump
through $Cl(\OO_x)$ via
$$\tk_{l+1,i} = \tk_{l,i}*\g^{s_{v(\tk_{l,i})}} \,, \quad \wk_{l+1,j} = \wk_{l,j}*\g^{s_{v(\wk_{l,j})}} \,,
\qquad (l \in \N_0 \,, \quad 1\leq i,\,j\leq m/2) \enspace .$$
The computations $\tk_{l,i} \to \tk_{l+1,i}$ and $\wk_{l,j}\to\wk_{l+1,j}$ are called {\em (kangaroo) jumps}. The {\em distance} of the $i$-th
tame kangaroo $T_i$ (or $j$-th wild kangaroo $W_j$) at step $l$ is the discrete logarithm of the ideal $\tk_{l,i}$ (or $\wk_{l,i}$) with
respect to the base ideal $\g$ and is denoted by $d_l(T_i)$ and $d_l(W_j)$, for tame and wild kangaroos, respectively. (We note that this
definition of distance is different from, and not to be confused with, the notion of infrastructure distance.) Specifically, we initialize
$d_0(T_i) = E+(i-1)\nu$ and $d_0(W_j) = (j-1)\nu$, for each $1\leq i,\,j \leq m/2$, so that
$$d_{l+1}(T_i) = d_l(T_i) + s_{v(\tk_{l,i})} \,, \quad d_{l+1}(W_j) = d_l(W_j) + s_{v(\wk_{l,j})} \,, \qquad (l \in \N_0) \enspace .$$
If $\tk_{A,i} = \wk_{B,j}$, for some $A,\,B\in\N_0$ and $1\leq i,\,j\leq m/2$, then we have a collision and $\g^{d_A(T_i)} = \g^{d_B(W_j)}$. If
$d_A(T_i) - d_B(W_j) \in [E-U,\,E+U]$, then we are guaranteed that $h = d_A(T_i) - d_B(W_j)$.

If there is a collision between any two kangaroos, then they will continue on the same path. Therefore, if there is a collision between two kangaroos
of the same herd, then we cannot obtain any information about $h$, so we must re-initialize one of the two kangaroos. Without loss of generality,
suppose that the two tame kangaroos $T_1$ and $T_2$ collide at the distance $d_A(T_1) = d_B(T_2)$. For a small $c \in \N$, set $\tk_{A+1,1} =
\tk_{A,1}*\g^c$ and $d_{A+1}(T_1) = d_A(T_1) + c$, then let $T_1$ continue jumping on its new path as usual. $T_2$ may continue along the same path
as before without interruption.
%In practice, there are few collisions, but a variation due to Pollard \cite{pollard} guarantees that no two
%kangaroos of the same herd collide. Theoretical analysis and implementation results of this method may be found in \cite{pollard} and \cite{st05}.

A key feature of the Kangaroo algorithm is that there is no need to store every jump. Using the idea of van Oorschot and Wiener \cite{vow}, we
only store {\em distinguished points}. In order to avoid confusion with the concept of distinguished divisors and ideals, such points will be
called {\em (kangaroo) traps} instead. To this end, define another hash function $z:\I(\OO_x) \to \{0,\,\ldots,\,\theta-1\}$. Install a trap, that
is, store a kangaroo $\kk$, if $z(\kk) = 0$. In this way, we expect to set a trap every $\theta$ jumps. If $\theta$ is sufficiently large, then the
storage requirement is very small.  Note that we only detect collisions between traps,  but since colliding kangaroos travel along the same
path following their first collision, a collision in a trap will eventually be found.

If it is known that there exist integers $a,\,b\in\N_0$ such that $0\leq a < b$ and $h\equiv a\ppmod{b}$, then we can make adjustments to the jump
set and initializations to only operate within the congruence class $a \ppmod{b}$. We change the estimate $E$ to $E - (E \ppmod{b}) + a$, so that
$E\equiv a\ppmod{b}$ for the revised value of $E$, and choose $\nu$ and the jump distances such that $b\mid \nu$ and $b\mid s_i$, for each $1\leq
i\leq 64$. The remaining initializations are the same.
%In this way, we have $d_l(T_i) \equiv a\ppmod{b}$ and $d_l(W_j) \equiv 0\ppmod{b}$, for any tame kangaroo $T_i$, with $1\leq i \leq m/2$, wild
%kangaroo $W_j$, with $1\leq j\leq m/2$, and $l\in\N_0$. Thus, if there is a collision, $\g^{d_A(T_i)} = \g^{d_B(W_j)}$, then $d_A(T_i) - d_B(W_j)
%\equiv a \ppmod{b}$.

In Algorithm \ref{alg:roo0}, we formalize the procedures described above.
%Its correctness follows from our discussion.
%Afterwards, we will optimize the expected running time of this algorithm, over all cubic function fields over $\Fq(x)$ and genus $g$, to
%justify the choices of particular variables.

\begin{algorithm}[ht]
\caption{Computing $h$ via the Kangaroo Algorithm - Signature $(3,1)$} \label{alg:roo0}
\begin{algorithmic}[1]
\REQUIRE A prime power $q$; monic, relatively prime, square-free $G,\, H \in \Fq[x]$ such that $3\nmid\deg(GH^2)$; $a,\,b \in \N_0$
such that $h \equiv a \ppmod{b}$ (or $b = 1$ and $a = 0$ if no non-trivial $b$ is known); integers $E,\,U\in\N$ such that $|h-E|<U$ and
$2(E-U) > E+U$; and an even integer $m$, the number of processors.%
\ENSURE The divisor class number $h$ of $K= \Fq(C)$, where $C: Y^3 = GH^2$.%
\STATE $g:= \deg(GH)-1$.%
\STATE Find an estimate, $\ah := \ah(q,\,g)$, of the expected value of $|h-E|/U$ via Table \ref{tab:alphas}.
\STATE $\beta := \left[(m/2)\sqrt{\ah b U}\right]$, $\nu := [2\beta/m] - \left([2\beta/m]\ppmod{b}\right)$.%
\STATE $\theta := 2^{\left[\lg(\beta)/2\right]}$, $E:= E - (E\ppmod{b}) + a$.%
\STATE Choose random integers $0<s_i\leq 2\beta$, for $1\leq i\leq 64$, with $Mean\left(\{s_i\}\right) = \beta$ and $b\mid s_i$.
\label{step:choose-jumps}%
\STATE Generate a random ideal $\g$.%
\STATE Define hash functions $v:\I \to \{1,\,\ldots,\,64\}$ and $z: \I \to \{0,\,\ldots,\,\theta-1\}$.%
\FOR {$i = 1$ to $m/2$}
   \STATE Initialize the tame kangaroo, $T_i$: $\tk_{0,i} := \g^{E+(i-1)\nu}$ and $d_0(T_i) := E+(i-1)\nu$.
   \STATE Initialize the wild kangaroo, $W_i$: $\wk_{0,i} := \g^{(i-1)\nu}$ and $d_0(W_i) := (i-1)\nu$.
\ENDFOR \IF{$z(\tk_{0,i}) = 0$ or $z(\wk_{0,i}) = 0$, for some $1\leq i \leq m/2$}
   \STATE Store the respective ideal and its distance.
\ENDIF%
\STATE $j:=0$.%
\WHILE{A collision between a tame and a wild kangaroo has not been found}
   \FOR{$i = 1$ to $m/2$}
      \STATE $\tk_{j+1,i} := \tk_{j,i}*\g^{s_{v(\tk_{j,i})}}$ and $d_{j+1}(T_i) := d_{j+1}(T_i) + s_{v(\tk_{j,i})}$.
      \STATE $\wk_{j+1,i} := \wk_{j,i}*\g^{s_{v(\wk_{j,i})}}$ and $d_{j+1}(W_i) := d_{j+1}(W_i) + s_{v(\wk_{j,i})}$.
   \ENDFOR
   \STATE $j:= j+1$.
   \IF{$z(\tk_{j,i}) = 0$ or $z(\wk_{j,i}) = 0$, for some $1\leq i \leq m/2$}
      \STATE Store the respective ideal and its distance.
   \ENDIF
\ENDWHILE
\IF{$\tk_{A,\,i} = \wk_{B,\,j}$}
   \RETURN $h := d_A(T_i)-d_B(W_j)$
\ENDIF
\end{algorithmic}
\end{algorithm}

The following analysis is a generalization of similar ideas in \cite{st02,st05}. It  justifies our choices of certain variables in Algorithm
\ref{alg:roo0}. Since it mainly depends on the set-up of the Kangaroo method and not on the underlying function field, the proof is omitted; for
details, see \cite{ericphd}. Henceforth, we denote by $[R]$ the nearest integer to a real number $R$.

\begin{proposition}
\label{prop:roo0} Let $K$ be a purely cubic function field of signature $(3,\,1)$. Suppose that there exist integers $a,\,b\in\N$ such that $h
\equiv a \ppmod{b}$. Then the expected heuristic running time, over all cubic function fields over $\Fq(x)$ of genus $g$, to compute $h$ via
Algorithm \ref{alg:roo0} is minimized by choosing an average jump distance of $\beta = \left[(m/2)\sqrt{\alpha b U}\right]$, where $m$ is the
(even) number of processors and $\alpha = \alpha(q,\,g) < 1/2$ is the mean value of $|h-E|/U$ over all cubic function fields $K$ over $\Fq(x)$
of genus $g$. For this choice of $\beta$, the total expected heuristic running time of Algorithm \ref{alg:roo0} for each kangaroo is
$(4/m)\sqrt{\alpha U/b} + \theta + O(1)$ ideal compositions as $q \to \infty$, where traps are set on average every $\theta$ jumps.
\end{proposition}

For further practical considerations, Stein and Teske \cite{teske, st05} note that for hyperelliptic function fields, if the jump distances
$s_1,\,\ldots,\,s_{64}$ are chosen randomly, then the number of useless collisions appears independent of the choice of the initial spacing
$\nu$, and suggest using $\nu \lessapprox 2\beta/m$. They also recommend choosing $s_i \leq 2\beta$ for $1\leq i\leq 64$, since such choices
yielded results which were slightly better than those using other upper bounds. We expect no difference for cubic function fields and therefore
chose $|J| = 64$ to be a power of $2$, so that the hash function $v$ is fast while still achieving a sufficient level of randomization, and
also the space to store the jumps is not too large. Lastly, we expect to store $O\left(\sqrt{U}/\theta\right)$ ideals.
Teske \cite{teske} suggests in the generic group setting taking $\theta = 2^{[\lg(\beta)/2]+c}$ for
some small integer $c$. For this choice, we have $\theta = O\left(\sqrt[4]{U}\right)$,
%which is not constant with respect to $\beta$, as assumed in the analysis of Proposition \ref{prop:roo0}. However,
%for the first two terms in the expression for $T_K(\beta)$ in \eqref{eq:roo0}, we have $\kh/\beta,\, (4\beta)/(bm^2) = O\left(\sqrt{U}\right)$,
%which dominate the term $\theta$. Therefore, this choice of $\theta$ does not significantly change the analysis in the proof of Proposition
%\ref{prop:roo0}. Since $T_K(\beta) = 4\sqrt{\alpha U/b} + \theta m + O(1) = O\left(\sqrt{U}\right)$, we expect to store
%$O\left(\sqrt{U}/\sqrt[4]{U}\right) = O\left(\sqrt[4]{U}\right)$ ideals,
which is a reasonable number in practice.
%We may choose a custom $\theta$ for each computation so that the number of ideals that we store
%is not too large, but so that the time between setting traps is not too long. Table 3 of \cite{st02} compares the results of experiments
%using varying $\theta$. For large examples, we used such values of $\theta$; specific choices are given in the examples in Section
%\ref{S:results}. %Based on these results, Tables \ref{tab:genus3est}, \ref{tab:genus4est}, \ref{tab:genus3estR}, and \ref{tab:genus4estR}
%give estimates for the expected number of traps for a certain choice of $\theta$ in cubic function fields over large base field, of genera
%$3$ and $4$, and signatures $(3,1)$ and $(1,1;1,2)$.

In the next section, we describe appropriate changes to use the Kangaroo method in the infrastructure of a cubic function field of signature
$(1,1;1,2)$. These changes correspond to the changes made for real hyperelliptic function fields.

\section{The Kangaroo Method in $\R$}
\label{S:rooinf}

If $K$ has signature $(1, 1; 1, 2)$, then we wish to determine $R_x$ via the computation of some multiple $h_0$ of $R_x$. Under the assumption
$2(E-U) > E+U$, if we find $h_0 \in [E-U,\,E+U]$, then in fact $h_0 = h$. In this case, we adapt the description of the Kangaroo algorithm in
Section \ref{S:kangaroo} to operate in $\R$ and show how to take advantage of the faster baby step operation. We formalize these modifications in
Algorithm \ref{alg:roo1} and note its running time in Proposition~\ref{prop:rooinf}.

To be consistent with earlier notation, a (tame or wild) kangaroo $Z$ in this context is a sequence of infrastructure divisors, and we write $Z
= \{\kk_0,\,\kk_1,\,\ldots\} \subseteq \R$. If $\kk_l \in Z$, then the distance of $Z$ at step $l$, $d_l(Z) = \delta(\kk_l)$, is the distance of
$\kk_l$ as defined for divisors in $\R$. We initialize each tame kangaroo $T_i$ ($1\leq i\leq m/2$) at the distinguished divisor $\tk_{0,i} =
D(2E + (i-1)\nu) \in \R$. Likewise, each wild kangaroo $W_j$ ($1\leq j\leq m/2$) is initialized at $\wk_{0,j} = D((j-1)\nu) \in \R$. In this
case, the jump set is $J = \{\g_1,\,\ldots,\,\g_{64}\}$, where $\g_i = D(s_i)$, for $1\leq i\leq 64$. However, since we do not necessarily have
$\delta(\g_i) = s_i$ for all $1\leq i\leq 64$, we store each distance, $\delta(\g_i)$, rather than the random integers $s_i$, for $1\leq i\leq
64$. Thus, for each step of the algorithm, we have $\kk_{l+1} = \kk_l \oplus \g_{v(\kk_l)}$, with the distances updated by $d_{l+1}(Z) = d_l(Z)
+ \delta\left(\g_{v(\kk_l)}\right) + \delta$, where $\delta = \deg(\psi)$ and $\psi$ is as given in \eqref{eq:gs}. In this adaptation, if a tame
kangaroo $T_i$ at step $A$ collides with a wild kangaroo $W_j$ at step $B$, then $\tk_{A,i} = \wk_{B,j}$. Thus, $\delta_A(T_i) \equiv
\delta_B(W_j) \ppmod{2R_x}$, so $h_0 = \left(\delta_i(T_A) - \delta_j(W_B)\right)/2$ is a multiple of $R_x$. We give details on how to
determine $R_x$ from $h_0$ in Section \ref{S:implementation4}.

In the infrastructure setting, however, we may take advantage of the fact that baby steps are faster than giant steps in $\R$ to speed up the
regulator computation by a factor of approximately $\sqrt{\tau_1/2}$, where $\tau_1 = T_G/T_B$ and $T_G$ and $T_B$ are the respective times to
compute a giant step and a baby step. The following is a slight change from the idea found in Section 4.1 of \cite{st02}. For a real number
$\tau\geq 1$, let $\SC_{\tau} \subseteq \R$ such that $|\R|/|\SC_{\tau}| \approx \tau$.\footnote{If $q$ is prime, then one possible choice for
$\SC_{\tau}$ is $\SC_{\tau} = \{D\in\R \st d(\Psi(-D))(0) < [q/\tau]\}$, where $d(\f) \in \Fq[x]$ is the denominator of the fractional ideal
$\f$.} After each kangaroo jump (a giant step), we take baby steps until a divisor in $\SC_{\tau}$ is found, then we make the next kangaroo
jump. Below, we outline the kangaroo algorithm and give specific choices for $\beta$ and $\tau$ to optimize its running time.

\begin{algorithm}[ht]
\caption{Computing $h$ via the Kangaroo Algorithm - Signature $(1,1; 1,2)$}
\label{alg:roo1}
\begin{algorithmic}[1]
\REQUIRE A prime power $q\equiv 2\ppmod{3}$; monic, relatively prime, square-free $G,\,H \in \Fq[x]$ such that $3\mid\deg(GH^2)$;
integers $E,\,U\in\N$ such that $|h-E|<U$ and $2(E-U) > E+U$, and an even integer $m$, the
number of processors.%
\ENSURE A multiple $h_0$ of the regulator $R_x$ of $K = \Fq(C)$, where $C: Y^3 = GH^2$.%
\STATE $g:= \deg(GH)-2$.%
\STATE Find an estimate, $\ah := \ah(q,\,g)$, of the expected value of $|h-E|/U$ via Table \ref{tab:alphas}.%
\STATE Determine the appropriate value of $\tau$ from Table \ref{tab:tau}.%
\STATE $\beta := \left[m\sqrt{(2\tau-1) \ah U}\right]-2(\tau-1)$, $\nu := [2\beta/m]$, $\theta := 2^{\left[\lg(\beta)/2\right]}$.
\STATE Choose random integers $g+2 \leq s_i\leq 2(\beta+\ph(g))+1$, with $1\leq i\leq 64$,
    such that $Mean\left(\{s_i\}\right) = \beta + \ph(g) + 1/2$.
\STATE Compute the jump set $J := \left\{D(s_1),\,\ldots,\, D(s_{64})\right\}$.% via Algorithm \ref{alg:below}.
\STATE Define hash functions $v:\I \to \{1,\,\ldots,\,64\}$ and $z: \I \to \{0,\,\ldots,\,\theta-1\}$.%
\FOR{$i = 1$ to $m/2$}
   \STATE Initialize the tame kangaroos, $T_i$: $\tk_{0,i} := D(2E+(i-1)\nu)$.
   \STATE Initialize the wild kangaroos, $W_i$: $\wk_{0,i} := D((i-1)\nu)$.
   \WHILE{$\tk_{0,i}\notin S_{\tau}$ }%$L(\Psi(\tk_{0,i}))(0) \geq [q/\tau]$}
      \STATE $\tk_{0,i} := bs(\tk_{0,i})$.% via Algorithm \ref{alg:babysteps}.
   \ENDWHILE
   \WHILE {$\wk_{0,i}\notin S_{\tau}$}%$L(\Psi(\wk_{0,i}))(0) \geq [q/\tau]$}
      \STATE $\wk_{0,i} := bs(\wk_{0,i})$
   \ENDWHILE
\ENDFOR

%\algstore{roobreak}

%\end{algorithmic}
%\end{algorithm}

%\begin{algorithm}[ht]
%\caption{Computing $h$ via the Kangaroo Algorithm - Signature $(1,1; 1,2)$, The Main Loop}
%\begin{algorithmic}[1]
%\REQUIRE A prime pwoer $q\equiv 2\ppmod{3}$; monic, relatively prime, square-free $G,\,H \in \Fq[x]$ such that $3\mid\deg(GH^2)$;
%$K = \Fq(C)$, $C: Y^3 = GH^2$,
%so that $\sig(K) = (1,1;1,2)$; integers $E,\,U\in\N$ such that $|h-E|<U$ and $2(E-U) > E+U$, and an even $m$, the number of processors.
%\ENSURE A multiple, $h_0$, of the regulator, $R_x$, of $K$.
%\algrestore{roobreak}

%\STATE $j:=k:=0$
\WHILE{A collision between a tame and a wild kangaroo has not been found}
   \FOR{$i = 1$ to $m/2$}
      \IF{$z(\tk_{0,i}) = 0$ or $z(\wk_{0,i}) = 0$}
         \STATE Store the respective divisor(s) and its (their) distance(s).
      \ENDIF
      \STATE $\tk_{j+1,i} := \tk_{j,i} \oplus D\left(s_{v(\tk_{j,i})}\right)$ and $\wk_{k+1,i} := \wk_{k,i} \oplus D
      \left(s_{v(\wk_{k,i})}\right)$ \label{step:wildgs} % via Algorithm~\ref{alg:giantsteps1}. \label{step:tamegs}
      %\STATE $\delta_{j+1}(T_i) := \delta(\tk_{j+1,i}) = \delta(\tk_{j,i}) + \delta_{T,j,i}$, where $\delta_{T,j,i}$ is the second output of
        %Algorithm \ref{alg:giantsteps1} with the input divisors given in Step \ref{step:tamegs}.
      %\STATE $\wk_{k+1,i} := \wk_{k,i} \oplus D\left(s_{v(\wk_{k,i})}\right)$ \label{step:wildgs}
      %\STATE $\delta_{k+1}(W_i) := \delta(\wk_{k,i}) = \delta(\wk_{k,i}) + \delta_{W,k,i}$, where $\delta_{W,k,i}$ is the second output of
        %Algorithm \ref{alg:giantsteps1} with the input divisors given in Step \ref{step:wildgs}.
   \ENDFOR
%   \STATE $j:= j+1$ and $k := k + 1$
   \FOR{$i = 1$ to $m/2$}
      \WHILE{$\tk_{j,i}\notin S_{\tau}$}
         \STATE $\tk_{j+1,i} := bs(\tk_{j,i})$%, $j:= j+1$
         %\STATE $\delta_{j+1}(T_i) := \delta(\tk_{j+1,i}) = \delta(\tk_{j,i}) + \delta_{T,j,i}$, where $\delta_{T,j,i}$ is the second output of
            %Algorithm \ref{alg:babysteps} with the input divisor $\tk_{j,i}$.
      \ENDWHILE
      %\WHILE{$L(\Psi(\wk_{k,i}))(0) \geq [q/\tau]$}
      \WHILE{$\wk_{0,i} \notin S_{\tau}$}
         \STATE $\wk_{k+1,i} := bs(\wk_{k,i})$%, $k:= k+1$
         %\STATE $\delta_{k+1}(W_i) := \delta(\wk_{k+1,i}) = \delta(\wk_{k,i}) + \delta_{W,k,i}$, where $\delta_{W,k,i}$ is the second output of
            %Algorithm \ref{alg:babysteps} with the input divisor $\wk_{k,i}$.
         %\STATE $k:= k+1$
      \ENDWHILE
   \ENDFOR
\ENDWHILE \IF{$\tk_{A,\,i} = \wk_{B,\,j}$}
   \RETURN $h_0 := (\delta_A(T_i)-\delta_B(W_j))/2$.
\ENDIF
\end{algorithmic}
\end{algorithm}

One key difference between the following result to optimize the expected running time of Algorithm \ref{alg:roo1} and Proposition
\ref{prop:roo0} is that we express the running time in terms of $T_G$, rather than in terms of the number of kangaroo jumps. The following is a
slight improvement of Equation (4.8) of \cite{st02} in the hyperelliptic case.

\begin{proposition}
\label{prop:rooinf} If $K/\Fq(x)$ is a purely cubic function field of signature $(1,1;1,2)$. Then assuming Heuristic \ref{the-heuristic}, the
expected heuristic running time, over all cubic function fields over $\Fq(x)$ of genus $g$, to compute a multiple $h_0$ of $R_x$ via Algorithm
\ref{alg:roo1} is minimized by choosing $\beta = \left[m\sqrt{(2\tau-1) \alpha U}\right] - 2(\tau-1)$, where $m$ is the (even) number of
processors, $\tau = T_G/T_B$, $T_G$ and $T_B$ are the respective times required to compute a giant step and a baby step in $\R$, and $\alpha =
\alpha(q,\,g) < 1/2$ is the mean value of $|h-E|/U$ over all cubic function fields over $\Fq(x)$ of genus $g$. With these choices, the expected
heuristic running time for each kangaroo is
$$\left(\frac{4}{m}\sqrt{\frac{\alpha U}{2\tau-1}} + \frac{\theta}{\tau} + O(1)\right)\left(2 - \frac{1}{\tau}\right)T_G \enspace ,$$
as $q \to \infty$, where traps are set on average every $\theta$ jumps.
\end{proposition}

As in the signature $(3,1)$ case, we discuss the reasons for the choice of certain other variables in Algorithm \ref{alg:roo1}. First, when choosing
values for the $s_i$, $1\leq i\leq 64$, there are a few considerations arising from the reduction required for giant steps, so that the average
jump distance is as close to $\beta$ as possible in practice. For each $s_i$, we have $\delta(D(s_i)) = s_i$ with probability roughly $1/2$ and
$\delta(D(s_i)) = s_i-1$ with probability roughly $1/2$, for sufficiently large $q$, by Heuristic \ref{the-heuristic}. Moreover, by Part 2 of
Heuristic \ref{the-heuristic}, we have $\left(d_{l-1}(Z) + s_{v(\kk_{l-1})}\right) - d_l(Z) = \ph(g)$ with probability $1 - O(1/q)$. Therefore,
to adjust for this ``headwind,'' as well as the average difference $s_i - \delta(D(s_i))$, we must choose the $s_i$ so that $(s_1 + \cdots +
s_{64})/64 = \beta + 1/2 + \ph(g)$. Likewise, we choose $s_i \leq 2\big(\beta+\ph(g)\big) +1$ so that each jump has distance bounded
above by $2\beta$ with probability close to $1$. Finally, we cannot have $0\in J$, otherwise a kangaroo will become permanently stuck at one
divisor if it hashes to $0$, so we must set a lower bound on the choices of the $s_i$ to avoid this situation. By Theorem 5.3.10 of
\cite{ericphd}, we have $1\leq \delta(bs(0))\leq g+2$, so $s = g+2$ is the smallest integer that guarantees that $D(s) \neq 0$. Therefore, we
must choose $s_i \geq g+2$ for all $1\leq i\leq 64$. With these choices of the $s_i$, the average jump distance is in practice as close to
$\beta$ as possible.

Table \ref{tab:tau} lists values of $\tau$ for various signature $(1,1;1,2)$ situations of genera $3\leq g\leq 7$. In each case, we computed the
ratios using $10^6$ baby steps and $10^6$ giant steps in a function field $\Fq(C)$ with $q = 10^8 + 7$ and $C:Y^3 = GH^2$ where $G$ and $H$ were
random distinct irreducible polynomials with $\deg(G) \geq \deg(H)$.
%As is the case in hyperelliptic function fields \cite{andreas01}, this experimental evidence indicates that the baby step
%operation is significantly faster than the giant step operation, and in most cases we have $[\tau] = g$.

%\begin{table}[ht]
%\caption{Giant Step to Baby Step Ratio in Signature $(1,1; 1,2)$ Infrastructure}
%$$\begin{array}{|c|c|c||c|}
%\hline
%g & \deg(G) & \deg(H) & \tau \\
%\hline
%2 & 2 & 2 & 2.96977 \\ % q = 10^8 + 7
%\hline
%3 & 4 & 1 & 2.92374 \\ %2.88781 \\ q = 10^8 + 7  1.57078\sqrt{U}
%\hline
%4 & 6 & 0 & 3.87316 \\ %3.88672 \\ q = 10^8 + 7  1.17516
%4 & 3 & 3 & 4.11812 \\ %4.07611 \\ q = 10^8 + 7  1.14465
%\hline
%5 & 5 & 2 & 5.29813 \\ %5.22904 \\ q = 10^8 + 7  1.02508
%\hline
%6 & 7 & 1 & 5.86166 \\ %5.82408 \\ q = 10^8 + 7  0.89324
%6 & 4 & 4 & 6.10144 \\ %6.21530 \\ q = 10^8 + 7  0.87708
%\hline
%7 & 9 & 0 & 7.50799 \\ %7.28357 \\ q = 10^8 + 7  0.70847
%7 & 6 & 3 & 7.72477 \\ %7.45021 \\ q = 10^8 + 7  0.69888
%\hline
%\end{array}$$
%\label{tab:tau}
%\end{table}

\begin{table}[ht]
\caption{Giant Step to Baby Step Ratio in $\R$}
$$\begin{array}{|c|c|c||c|||c|c||c|}
\hline
g & \deg(G) & \deg(H) & \tau & \deg(G) & \deg(H) & \tau\\
\hline
2 & 2 & 2 & 2.96977 & & & \\
\hline
3 & 4 & 1 & 2.92374 & & & \\
\hline
4 & 6 & 0 & 3.87316 & 3 & 3 & 4.11812 \\
\hline
5 & 5 & 2 & 5.29813 & & & \\
\hline
6 & 7 & 1 & 5.86166 & 4 & 4 & 6.10144 \\
\hline
7 & 9 & 0 & 7.50799 & 6 & 3 & 7.72477\\
\hline
\end{array}$$
\label{tab:tau}
\end{table}

In the next section, we review the method of \cite{ss07} implemented here to compute the divisor class number of a cubic function field.

\section{Approximating $h$}
\label{S:EU}

\subsection{Idea of the Algorithm}
\label{S:idea}

Algorithm \ref{alg:h-idea} lists the three main phases of the method of \cite{ss07} to compute the divisor class number $h$ of a cubic function
field, with a fourth step in the signature $(1,1;1,2)$ case if the regulator $R_x$ is desired.
\begin{algorithm}[ht]
\caption{Computing $h$ and/or $R_x$ - The Idea} \label{alg:h-idea}
\begin{algorithmic}[1]
\STATE Compute an estimate $E$ of $h$ and an upper bound $U$ on the error $|h - E|$ so that $h \in [E-U,\,E+U]$.%
\STATE Determine extra information about $h$ such as congruences or divisibility by small primes or distribution of $h$ in the interval
$[E-U,\, E+U]$.%
\STATE Use the Baby Step-Giant Step method or the Kangaroo method to find $h$ in $[E-U,\,E+U]$ using $O\left(\sqrt{U}\right)$ ideal
compositions.%
\STATE If $K$ has signature $(1,1;1,2)$, then factor $h$ and let $R_x$ be the smallest factor $R'$ of $h$ such that $D(2R') = 0$.
\end{algorithmic}
\end{algorithm}

\subsection{Results and Notation for Phase 1}

For full details on the derivation of $E$ and $U$, we refer to \cite{sw,st02,st05}, \cite{ss07}, and \cite{ss08} for the cases of quadratic,
cubic, and arbitrary function fields, respectively. The idea, however, is to write $h$ as an infinite product over the places of $K$ via the
zeta function of $K$; $E$ is determined by computing the product up to a certain degree bound $\lambda$, and $U$ is determined by setting an
upper bound on the size of the tail.

Following \cite{ss07,ss08}, let
$$(x_1, \, x_2) = \left\{
\begin{array}{ll}
  (0, \,0)   & \mbox{if } \sig(K) = (3,1)\enspace ,\\% \infty=\infty_0^3 \enspace ,\\
  %(\iota, \,\iota^2)   & \mbox{if }\infty=\infty_0 \enspace ,\\
  (1, \,-1)  & \mbox{if } \sig(K) = (1,1;1,2)\enspace ,\\%\infty=\infty_0\infty_1 \enspace ,\\
  %(1, \,0) & \mbox{if } \sig(K) = (1,1;2,1)\enspace ,\\%\infty=\infty_0\infty_1^2 \enspace ,\\
  %(1, \,1) & \mbox{if } \sig(K) = (1,1;1,1;1,1)\enspace ,\\%\infty=\infty_0\infty_1\infty_2 \enspace .
\end{array}
\right.$$
and $A(K) = (g+2)\log(q) - \log((q -x_1)(q - x_2))$. Next, let $\iota$ be a primitive cube root of unity in some algebraic closure of
$\Fq$. If $P \in \Fq[x]$ is a monic irreducible polynomial, then let
\begin{equation}
\label{eq:Zfinite}
(z_1(P), \, z_2(P)) = \left\{ \begin{array}{ll}
(0,\,0)     & \mbox{if } \ideal{P}=\p^3 \enspace ,\\
(1, \,-1)   & \mbox{if } \ideal{P}=\p\q \enspace ,\\
(1, \,1)    & \mbox{if } \ideal{P}=\p\p'\p''\enspace ,\\
(\iota,\,\iota^2) & \mbox{if } \ideal{P}=\p \enspace ,
\end{array}
\right.
\end{equation}
and
\begin{equation} \label{eq:Snu}
S_{\nu}(n) = \sum_{\deg(P) = \nu}\left(z_1^n(P) + z_2^n(P)\right) \enspace .
\end{equation}
From (4.11) and (4.12) of \cite{ss07}, we have
\begin{equation}
\label{eq:h2}
h = \frac{q^{g+2}}{(q-x_1)(q-x_2)}\prod_{\nu=1}^{\infty}\prod_{\deg(P) = \nu}\frac{q^{2\nu}}{(q^{\nu} - z_1(P))(q^{\nu}-z_2(P))} \enspace .
\end{equation}
Equation \eqref{eq:h2} gives rise to the following result (see \cite{ss07}):

\begin{theorem}
\label{thm:logh} Let $K/\Fq(x)$ be a purely cubic function field of genus $g$ such that $\ch(K) \neq 3$. Then
$$\log(h) = A(K) + \sum _{n=1}^{\infty}\frac{1}{nq^{n}}\sum_{\nu\mid n}\nu S_{\nu}\left(\frac{n}{\nu}\right) \enspace .$$
\end{theorem}

\subsection{Determining $E$ and $U$}
\label{S:estimates}

We review how to explicitly calculate two estimates $E_1, E_2$ of $h$, and three bounds $U_1, U_2, U_3$ on the error, $|h-E|$,
as found in \cite{ss07, ss08}. For $i=1,\,2,\,3$ and a fixed $\lambda\in\N$, we write $\log(h) = \log(E_i'(\lambda)) + B_i(\lambda)$ so that $h
= E_i'(\lambda)e^{B_i(\lambda)}$. (For ease of notation, we write $E_2(\lambda)=E_3(\lambda)$ and $B_2(\lambda)=B_3(\lambda)$.) We then find a
sharp upper bound $\psi_i(\lambda)$ on $|B_i(\lambda)|$ and define $E_i(\lambda) := \left[E_i'(\lambda)\right]$ and $U_i(\lambda) := \left[
E'_i(\lambda)\left(e^{\psi_i(\lambda)} - 1\right)\right]$, so that $|h - E_i(\lambda)| \leq U_i(\lambda)$.

First, we have
$$\log(E_1'(\lambda)) = A(K) + \sum_{n=1}^{\lambda}\frac{1}{nq^{n}}\sum_{\nu\mid n}\nu S_{\nu}\left(\frac{n}{\nu}\right)\,, \quad
B_1(\lambda) = \sum_{n=\lambda + 1}^{\infty}\frac{1}{nq^{n}}\sum_{\nu\mid n}\nu S_{\nu}\left(\frac{n}{\nu}\right) \enspace .$$
Then
$$\psi_1(\lambda) = 2g\left(\log\left(\frac{\sqrt{q}}{\sqrt{q} - 1}\right) - \sum_{n=1}^{\lambda}\frac{1}{nq^{n/2}}\right) +
2\log\left(\frac{q}{q-1}\right) - 2\sum_{n=1}^{\lambda}\frac{1}{nq^n} $$
is a sharp upper bound on $|B_1(\lambda)|$. By moving some terms from $B_1(\lambda)$ to $E_1(\lambda)$, we obtain the following second estimate
$E_2(\lambda)$ and error bound. $E_2(\lambda)$:
$$\log(E_2'(\lambda)) =  A(K) +
\sum_{n=1}^{\infty}\frac{1}{nq^n}\sum_{\substack{\nu\mid n\\ \nu \leq \lambda}}\nu S_{\nu}\left(\frac{n}{\nu}\right) \,, \quad
B_2(\lambda) = \sum_{n=\lambda+1}^{\infty}\frac{1}{nq^n}\sum_{\substack{\nu\mid n \\ \nu>\lambda}}\nu S_{\nu}\left(\frac{n}{\nu}\right)
\enspace .$$
A sharp upper bound $\psi_2(\lambda)$ of $|B_2(\lambda)|$ is then given by
\begin{align*}
\psi_2(\lambda) &= \frac{2}{(\lambda + 1)}\left(gq^{-(\lambda+1)/2} + q^{-(\lambda+1)}\right) + \frac{2q}{(q-1)(\lambda+1)}\,
q^{-(\lambda+1)}\left(q^{(\lambda+1)/l} - 1\right)\notag\\
&+ \frac{2g}{(\lambda+2)}\frac{\sqrt{q}}{\left(\sqrt{q} - 1\right)}\,q^{-(\lambda+2)/2} + \frac{4}{(\lambda+2)}\frac{q}{(q-1)}
\frac{q^{(l-1)/l}}{(q^{(l-1)/l}-1)}\,q^{-(\lambda+2)(l-1)/l}  \enspace ,
\end{align*}
where $l$ is the smallest prime factor of $\lambda+1$.

Finally, we use extra information to obtain a sharper bound $\psi_3(\lambda)$ on $B_2(\lambda)$. Specifically, we can easily calculate $\nu
S_{\nu}((\lambda+1)/\nu)$ for all $\nu\mid(\lambda+1)$ such that $\nu\neq\lambda+1$. Thus, we obtain
\begin{align} \label{eq:B3}
\psi_3(\lambda) &= \frac{2g}{(\lambda+1)}\,q^{-(\lambda+1)/2}
   + \frac{q^{-(\lambda+1)}}{(\lambda + 1)} \left(2 + \left| \sum_{\substack{\nu\mid(\lambda+1)\\ \nu \neq \lambda+1}}
\nu S_{\nu}\left(\frac{\lambda+1}{\nu}\right)\right|\,\right) \\
& \qquad + \frac{2g}{(\lambda+2)}\frac{\sqrt{q}}{\left(\sqrt{q} - 1\right)}\,q^{-(\lambda+2)/2}  \nonumber \\
& \qquad + \frac{4}{(\lambda+2)} \frac{q}{(q-1)}\frac{q^{(l-1)/l}}{(q^{(l-1)/l}-1)}\,q^{-(\lambda+2)(l-1)/l}  \nonumber \enspace .
\end{align}

\subsection{Complexity and Optimization}
\label{S:runtime} For complete details on the analysis of the running time of Algorithm \ref{alg:h-idea}, we refer to \cite{ss07, ss08}. Here,
we simply state that as $q\to\infty$, the optimal choice for $\lambda$ is
\begin{equation}
\label{eq:lambda}
\lambda = \left\{\begin{array}{ll}
\lfloor (2g-1)/5 \rfloor & \mbox{if } g \equiv 2 \ppmod{5} \enspace , \\
\left[(2g-1)/5\right]    & \mbox{otherwise} \enspace .
\end{array}\right.
\end{equation}
If $g \leq 2$, then $\lambda = 0$, so the estimate in Phase 1 is completely determined by the infinite component and runs in polynomial time; there
is no asymptotic improvement over using the Hasse-Weil bounds. However, if $g \geq 3$, then we have the following result.

\begin{theorem}
\label{thm:h-idea} If $K$ is a cubic function field of genus $g\geq 3$, then the complexity of Algorithm \ref{alg:h-idea} is
$O\left(q^{[(2g-1)/5]+\varepsilon(g)}\right)$ ideal or infrastructure compositions, as $q \to \infty$, where
$$\varepsilon(g) = \left\{\begin{array}{rl}
0    & \mbox{if } g \equiv 0,\,3 \ppmod{5} \enspace ,\\
1/4  & \mbox{if } g \equiv 1     \ppmod{5} \enspace ,\\
-1/4 & \mbox{if } g \equiv 2     \ppmod{5} \enspace ,\\
1/2  & \mbox{if } g \equiv 4     \ppmod{5} \enspace .\\
\end{array}\right.$$
\end{theorem}

In practice, Step 2 of Algorithm \ref{alg:h-idea} requires a negligible amount of time since the information associated with this step is known
in advance. Also, Step 4 is faster than Steps 1 and 3 since factoring is asymptotically faster than the overall running time of the
algorithm. Therefore, the overall complexity of Algorithm \ref{alg:h-idea} is found by balancing the running times of Steps 1 and 3. However,
Step 1 requires polynomial arithmetic at each step, whereas Step 3 requires ideal or infrastructure arithmetic at each step, which is much
slower. As a result, Step 3 dominates the overall running time in practice.

Next, we discuss practical issues surrounding actual implementations of each step of Algorithm \ref{alg:h-idea}. We remark that this is the
first time that this algorithm has been implemented for cubic function fields.

\section{Implementation Details}
\label{S:details}

\subsection{Implementation Details for Phase 1}
\label{S:implementation1} This section presents a number of algorithms and results to apply to the problem of computing both approximations
$E_1$ and $E_2$ of $h$ in Step 1 of Algorithm \ref{alg:h-idea}. All the algorithms and derivations in this section are new.
Henceforth, we assume that $q$ is prime. While the following methods can be extended to composite $q$, we make this restriction to facilitate
the execution of Phase 1 since there are very straightforward ways to loop through the set of all irreducible polynomials up to a fixed degree.
In addition, ideal, infrastructure, and polynomial arithmetic is faster using $q$ prime .

First, we note that the splitting behavior of an ideal $\ideal{P}$ in $\OO_x$, where $P \in \Fq[x]$ is irreducible), can be ascertained by
computing the cubic power residue symbol $[GH^2/P]_3$ if $q \equiv 1\ppmod{3}$. This can be done via Algorithm 6.2 of \cite{ss07} and has
essentially the same complexity as the Euclidean Algorithm applied to $GH^2$ and $P$. If $q \equiv 2\ppmod{3}$ and $\deg(P)$ is even, we have
$[P/Q]_3 \equiv Q^{\left(q^{\deg(P)}-1\right)/3} \ppmod{P}$, which is computed using $O(\deg(P)\log(q))$ polynomial operations.

The following equation is the core of Step 1 of Algorithm \ref{alg:h-idea}, and computes $z_1(P)^n + z_2(P)^n$, where $P\in\Fq[x]$ is
irreducible, $n\in\N$, and $z_1(P)$ and $z_2(P)$ are defined as in \eqref{eq:Zfinite}. Combining \eqref{eq:Zfinite} with Theorem
\ref{thm:primes}, and setting $\chi(P) = [GH^2/P]_3$, we have:
\begin{equation}
\label{eq:zs}
z_1(P)^n + z_2(P)^n = \left\{\begin{array}{rl}
-1 & \mbox{ if } q^{\deg(P)}\equiv 1\ppmod{3},\, \chi(P)\neq 1,\mbox{ and } 3\nmid n \,, \\
 0 & \mbox{ if } P\mid GH \mbox{ or } q^{\deg(P)}\equiv 2\ppmod{3} \mbox{ and } 2\nmid n \,, \\
 2 & \mbox{ otherwise} \,.
\end{array}\right.
\end{equation}

Next, we compute the value of $S_{\nu}(n)$ as given in \eqref{eq:Snu} for $1 \leq \nu \leq \lambda$ and all $n\in \N$. From
\eqref{eq:Zfinite} and \eqref{eq:zs}, we see that $S_{\nu}(n) = S_{\nu}(6k+n)$ for all $k,\,n\in\N$. Moreover, if $q^{\deg(P)}\equiv
1\ppmod{3}$, then $S_{\nu}(1) = S_{\nu}(2) = S_{\nu}(4) = S_{\nu}(5)$ and $S_{\nu}(3) = S_{\nu}(6)$; and if $q^{\deg(P)}\equiv 2\ppmod{3}$, then
$S_{\nu}(1) = S_{\nu}(3) = S_{\nu}(5) = 0$ and $S_{\nu}(2) = S_{\nu}(4) = S_{\nu}(6)$. For the cases $q^{\deg(P)} \equiv 2\ppmod{3}$ and $2\mid n$, or $q^{\deg(P)} \equiv 1\ppmod{3}$ and $3\mid n$, we have $z_1(P)^n + z_1(P)^n = 0$ if $P|GH$ and $z_1(P)^n + z_1(P)^n = 2$ otherwise.
In light of this, let $I_{\nu}$ be the number of irreducible polynomials of degree $\nu$ and $F_{\nu}$ the number of prime divisors of $GH^2$ of
degree $\nu$. Using well-known formulas for $I_{\nu}$, we have
\begin{align}
\label{eq:Snneq1}
S_{\nu}(n) = 2(I_{\nu}-F_{\nu}) = 2\left(\frac{1}{\nu}\sum_{d\mid \nu}\mu\left(\frac{\nu}{d}\right)q^d - F_{\nu}\right)%\notag\\
%= 2\left(\frac{1}{\nu}\left(q^{\nu} + \sum_{\substack{d\mid \nu\\d\neq \nu}}\mu\left(\frac{\nu}{d}\right)q^d\right) - F_{\nu}\right) \enspace ,
\end{align}
for these cases, where $\mu$ is the M\"{o}bius function. If $q^{\deg(P)} \equiv 1 \ppmod{3}$, then we must compute $S_{\nu}(1)$ by determining
the splitting behavior of each irreducible polynomial of degree $\nu$. Algorithms \ref{alg:computeS1} and \ref{alg:computeS2} check each
irreducible polynomial for the cases $\deg(P)= 1$ and $\deg(P)=2$, respectively. We also note that Algorithms \ref{alg:computeS1} and
\ref{alg:computeS2} may be parallelized by letting each processor run on distinct blocks of the interval $0\leq c < q$.

\begin{algorithm}[ht]
\caption{Computing $S_1(1)$} \label{alg:computeS1}
\begin{algorithmic}[1]
\REQUIRE A prime $q$, and monic, relatively prime, square-free $G,\,H\in\Fq[x]$.%
\ENSURE $S_1(1)$.%
\STATE $P := x$, $c := 0$, $S := 0$.%
\WHILE {$c < q$}
   \STATE Determine $z:= z_1(P)^n+z_2(P)^n$ as in \eqref{eq:zs}.
   \STATE $S := S + z$, $c := c+1$, $P := x + c$.
\ENDWHILE%
\RETURN $S_1(1) = S$.
\end{algorithmic}
\end{algorithm}

\begin{algorithm}[ht]
\caption{Computing $S_2(1)$} \label{alg:computeS2}
\begin{algorithmic}[1]
\REQUIRE A prime $q$, and monic, relatively prime, square-free $G,\,H\in\Fq[x]$.%
\ENSURE $S_2(1)$.%
\STATE $c := 1$, $j := 0$, $S := 0$.%
\WHILE {$c < q$}
   \IF {the Legendre symbol $(c/q) \neq 1$ } \label{step:nonsquare}
      \WHILE {$j < q$}
         \STATE $P:=(x-j)^2+c$.
         \STATE Determine $z:= z_1(P)^n+z_2(P)^n$ as in \eqref{eq:zs}.
         \STATE $S := S + z$, $j := j+1$.
      \ENDWHILE
   \ENDIF
   \STATE $j := 0$, $c := c+1$.
\ENDWHILE \RETURN $S_2(1) = S$.
\end{algorithmic}
\end{algorithm}

Finally, we give equations to determine $E_1$ and $E_2$. To compute $E_1$, we simply evaluate the sum given in Section \ref{S:estimates}. To
compute $E_2$, we evaluate the sum in Section \ref{S:estimates} by reversing the order of summation. Let $\nu m = n$ so that
$$\log{E_2'(\lambda)} = A(K) + \sum_{\nu=1}^{\lambda}\sum_{m=1}^{\infty}\frac{\nu S_{\nu}(m)}{\nu mq^{\nu m}} \enspace .$$
Using the identity $\sum_{m \in \mathbb{N}} 1/kmq^{km} = (1/k)\log\left(q^k/(q^k-1)\right)$, we have
\begin{align}
\label{eq:E2a}
\log(E_2') = A(K) & + \sum_{\nu=1}^{\lambda}\left(-S_{\nu}(1)\log\left(\frac{q^{\nu}-1}{q^{\nu}}\right)\right.\notag \\
&\quad+ \left.\frac{1}{3} (S_{\nu}(1)-S_{\nu}(3))\log\left(\frac{q^{3\nu}-1}{q^{3\nu}}\right)\right) %\enspace ,
\end{align}
if $q \equiv 1\ppmod{3}$. If $q \equiv 2\ppmod{3}$, then %$\log(E_2'(\lambda))$ is
\begin{align}
\label{eq:E2a2}
\log(E_2'(\lambda)) = A(K) & + \sum_{m=1}^{\lfloor\lambda/2\rfloor}\left(-S_{2m}(1)\log\left(\frac{q^{2m}-1}{q^{2m}}\right)\right.\notag \\
&\quad+  \left.\frac{1}{3}(S_{2m}(1)-S_{2m}(3))\log\left(\frac{q^{6m}-1}{q^{6m}}\right)\right)\notag \\
&\quad+ \sum_{m=1}^{\lfloor(\lambda+1)/2\rfloor}\left(-S_{2m-1}(1)\log\left(\frac{q^{2m-1}-1}{q^{2m-1}}\right)\right. \notag \\
&\quad+ \left.\frac{1}{2}(S_{2m-1}(1)-S_{2m-1}(2))\log\left(\frac{q^{4m-2}-1}{q^{4m-2}}\right)\right) \enspace .
\end{align}

\subsection{Implementation Details for Phase 2}
\label{S:implementation2}

For Phase 2 of Algorithm \ref{alg:h-idea}, we gather extra information about $h$ to effectively reduce the size of the interval, $[E-U,\,E+U]$,
determined in Phase 1. One observation is that $h$ is not uniformly distributed in this interval, and tends to be close to the approximation
$E$. In Sections \ref{S:kangaroo} and \ref{S:rooinf}, we described how to apply the average $\alpha(q,\,g) = Mean(|h-E|/U)$, taken over all
cubic function fields over $\Fq(x)$ of genus $g$, to optimize the expected running time of the  Kangaroo algorithms. Specifically the Kangaroo
algorithm is optimized with a shorter average jump length in order to concentrate our effort on the middle of the interval $[E-U,\,E+U]$,
thereby obtaining a speed-up by a factor of $(1+2\alpha(q,\,g))/\left(2\sqrt{2\alpha(q,\,g)}\right)$.

However, values of $\alpha(q,\,g)$ are very difficult to compute precisely, so in practice we apply an approximation $\ah(q,\,g)$ of $\alpha(q,\,g)$
instead. Table \ref{tab:alphas} in Section \ref{S:results} lists approximations $\ah(q,\,g)$ for selected values of $q$ and $g$, based on a
large sampling of cubic function fields of characteristic $q$ and genus $g$. For a fixed genus $g$, we assume that there is a limiting value
$\alpha(g) = \lim_{q\to\infty}\alpha(q,\,g)$, as is the case for hyperelliptic function fields \cite{st02b}, so that in practice, we can
interpolate or extrapolate as needed when applying these approximations for a given $q$ in Phase 2. As such, the information for this phase is
determined in advance. We will discuss the values of $\alpha(q,\,g)$ in more depth in Section \ref{S:results}.

\subsection{Implementation Details for Phase 3}
\label{S:implementation3}

Algorithms \ref{alg:roo0} and \ref{alg:roo1} presented earlier implement Phase 3 of Algorithm \ref{alg:h-idea}.

A more detailed description of the outline given in Algorithm \ref{alg:h-idea} for computing the divisor class number of a purely cubic function
field is provided in Algorithm \ref{alg:classnumber} below.

\begin{algorithm}[ht]
\caption {Class Number Computation for Purely Cubic Function Fields} \label{alg:classnumber}
\begin{algorithmic}[1]
\REQUIRE A prime $q$; monic, relatively prime, square-free $G,\, H \in \Fq[x]$; and $K = \Fq(C)$, where $C: Y^3 = GH^2$.%
\ENSURE The divisor class number $h$ of $K$.%
\IF{$3\mid\deg(GH^2)$}
   \STATE $g:=\deg(GH)-2$
\ELSE
   \STATE $g:=\deg(GH)-1$.
\ENDIF%
\STATE Set $\lambda$ via \eqref{eq:lambda}.%
\FOR {$\nu = 1$ to $\lambda$}
   \IF{$q^{\nu} \equiv 1\ppmod{3}$}
      \STATE Compute $S_{\nu}(1)$ via Algorithm \ref{alg:computeS1}, \ref{alg:computeS2}, etc.
      \STATE Compute $S_{\nu}(3)$ via \eqref{eq:Snneq1}.
   \ELSE
      \STATE $S_{\nu}(1) := 0$
      \STATE Compute $S_{\nu}(2)$ via \eqref{eq:Snneq1}.
   \ENDIF
\ENDFOR%
\IF {$q\equiv 1\ppmod{3}$}
   \STATE Compute $E := [\exp(\log(E_2'))]$ via \eqref{eq:E2a}, $\psi_3$ via \eqref{eq:B3}, and $U := \left[E_2'\left(e^{\psi_3}-1\right)\right]$.
   \STATE Compute and output $h$ via Algorithm \ref{alg:roo0}.
\ELSE
   \STATE Compute $E := [\exp(\log(E_2'))]$ via \eqref{eq:E2a2}, $\psi_3$ via \eqref{eq:B3}, and $U := \left[E_2'\left(e^{\psi_3}-1\right)\right]$.
   \STATE Compute and output $h$ via Algorithm \ref{alg:roo1}.
\ENDIF
\end{algorithmic}
\end{algorithm}

Theorem \ref{thm:h-idea} implies the following.

\begin{theorem}
With $-1/4 \leq \varepsilon(g) \leq 1/2$ as in Theorem \ref{thm:h-idea}, the complexity of Algorithm \ref{alg:classnumber} is $O\left(q^{[(2g-1)/5]
+\varepsilon(g)}\right)$ ideal operations.
\end{theorem}

\subsection{Implementation Details for Phase 4}
\label{S:implementation4}

Algorithm \ref{alg:S-regulator1} outlines the procedure for the final phase of Algorithm \ref{alg:h-idea} for the signature $(1,1;1,2)$ case; that
is, determining the regulator $R_x$, given a multiple $h_0$. We follow the procedure described in Algorithm 4.4 of \cite{sw}, making adaptations to
the cubic function field case. This technique uses the fact that the regulator $R_x$ is the smallest factor of $h_0$ such that
$D(2R_x) = 0$. The algorithm is an infrastructure analogue to determining the order of a group element from the group order.

\begin{algorithm}[ht]
\caption{Computing the Regulator of a Purely Cubic Function Field of Signature $(1,1; 1,2)$: Phase 4} \label{alg:S-regulator1}
\begin{algorithmic}[1]
\REQUIRE A multiple $h_0$ of $R_x$, a lower bound $l$ of $R_x$, a prime $q$, and monic, relatively prime, square-free $G,\, H \in \Fq[x]$.
\ENSURE The regulator $R_x$ of $K = \Fq(C)$, where $C: Y^3 = GH^2$.%
\STATE $h^* := 1$.%
\STATE Factor $h_0 = \prod_{i=1}^k p_i^{a_i}$. \label{step:factor}%
\FOR {$i=1$ to $k$}
   \IF {$p_i < h_0/l$}
      \STATE Find $1\leq e_i \leq a_i$ minimal such that $D\left(2h_0/p_i^{e_i}\right) \neq 0$. \label{step:below}% via Algorithm \ref{alg:below}.
      \STATE $h^* := p_i^{e_i-1}h^*$
   \ENDIF
\ENDFOR
\RETURN $R_x := h_0/h^*$
\end{algorithmic}
\end{algorithm}

We briefly comment on the running time of Algorithm \ref{alg:S-regulator1} relative to the running time of Algorithm~\ref{alg:classnumber},
especially in light of the factorization in Step 2. First, current heuristic methods to factor the integer $h_0$ require a subexponential
number of bit operations in $\log(h_0)$ using the Elliptic Curve Method \cite{ecm}, the Quadratic Sieve \cite{qs,mpqs,siqs}, or the General
Number Field Sieve \cite{gnfs} to achieve this running time. Furthermore, the loop in Steps 3-7 only requires a polynomial number (in $g$ and
$\log(q)$) of infrastructure operations. Therefore, determining $R_x$ from $h_0$ does not dominate the overall running time of Algorithm
\ref{alg:h-idea}. The largest divisor class numbers that we found have $28$ digits, which required only a few seconds to factor. In fact, we
simply used a basic implementation of Pollard's Rho factoring method~\cite{pollardrho}.

\section{Computational Results}
\label{S:results}

In this section, we present results and data obtained from the implementation of Algorithms \ref{alg:classnumber} and \ref{alg:S-regulator1} on
cubic function fields of signatures $(3,1)$ and $(1,1;1,2)$ and genera $4 \leq g \leq 7$. We first give experimental results that allowed us to
obtain constant-time speed-ups of Algorithm \ref{alg:h-idea}. We then discuss the problem of computing $\alpha_i(q,\,g) = Mean(|h-E_i|/U_i)$, where
the average is considered over all cubic function fields over $\Fq(x)$ of genus $g$. Finally, we list results of divisor class number and
regulator computations. For timing and technical considerations, we implemented our algorithms in C++ using NTL, written by Shoup \cite{ntl},
compiled using {\tt gcc}, and run on Sun workstations with AMD Opteron 148 $2.2$ GHz processors and $1$ GB of RAM running Fedora 7 Linux.

\subsection{General Optimization Data}
\label{S:optimize}

For this section, we applied the Baby Step-Giant Step method to $10,000$ function fields of signature $(3,1)$ of a fixed characteristic $q$ and
genus $g$, and organized the data from these computations to optimize implementations of Algorithm \ref{alg:classnumber}. This data provides
means to obtain a constant-time improvement over more straightforward implementations of this algorithm. First, we compared the accuracy of the
estimates $E_1$ and $E_2$. We then considered the minimal and maximal values of $|h-E_i|/U_i$ for each $i=1,\,2,\,3$, $q$, and $g$ to provide
further analysis of the estimates and compare the sharpness of the bounds $U_i$. Finally, for selected $q$ and for genera $3\leq g\leq 7$, we
list approximations $\ah_i(q,\,g)$ of $\alpha_i(q,\,g)$.

In each table of this and later sections, $\lambda$ is the degree bound used to compute the estimates $E_1$ and $E_2$, and $n$ is the number of
randomly chosen fields $K$ of the given characteristic $q$ and genus $g$ that we used in each experiment. In Table \ref{tab:estimates}, we
compare how well the two estimates $E_1$ and $E_2$ approximate $h$. Here, $\pm gs$ gives the average difference between the respective number of
giant steps computed using estimates $E_1$ and $E_2$, $\pm gs\%$ is the average percentage of the giant step
time gained or lost by using $E_2$ versus $E_1$, and $P_2$ is the percentage of the trials in which $E_2$ was the better estimate.

\begin{table}[ht]
\caption{Comparison of the Estimates $E_1$ and $E_2$}
$$\begin{array}{|rrr||r|r||c||c|}
\hline
q & g & \lambda  & \pm gs & \pm gs\% & P_2 &  n \\
\hline
997    & 3 & 1 &    0.23 &  0.054\% & 51.17\% & 10000\\
10009  & 3 & 1 &    3.24 &  0.061\% & 51.52\% & 10000\\
100003 & 3 & 1 &    7.52 &  0.014\% & 52.04\% & 10000\\
\hline
997    & 4 & 1 &   -9.52 & -0.067\% & 49.53\% & 10000\\
10009  & 4 & 1 &  140.01 &  0.033\% & 49.96\% & 10000\\
\hline
97     & 5 & 2 &   35.96 &  3.617\% & 53.86\% & 10000\\
997    & 5 & 2 & -106.92 & -0.181\% & 49.40\% & 10000\\
\hline
97     & 6 & 2 &  136.51 &  1.280\% & 52.23\% & 10000\\
463    & 6 & 2 &  530.11 &  0.139\% & 50.56\% & 10000\\
\hline
19     & 7 & 2 &   68.24 &  5.376\% & 53.65\% & 10000\\
97     & 7 & 2 & 1067.35 &  0.913\% & 51.27\% & 10000\\
\hline
\end{array}$$
\label{tab:estimates}
\end{table}

%We expected the difference between the two estimates to be more pronounced, but in each case, there was not much difference between the two. In all
%but the cases $g=4$ and $q = 997$, and $g=5$ and $q=997$, the second estimate, $gs$ was positive, and in all but the two genus $4$ experiments and
%$g=5$ and $q=997$, $E_2$ was the better estimate in most trials. We believe that the reason for some of the averages favoring $E_1$, while most
%favor $E_2$, is due to statistical variation in any sampling and the fact that there is not a very significant advantage to using $E_2$ over $E_1$.
%For larger examples, using either estimate would not lead to vastly different running times, as the $\pm gs$ and $\pm gs\%$ columns seem to indicate.
%For the genus $5$, $6$, and $7$ examples, the data suggest that the higher percentages in the $\pm gs\%$ column are due to the small base fields.
%Overall, $E_2$ is the better estimate, as expected, though not by a large margin.

In Table \ref{tab:minmaxalphas}, we give the minimum and maximum values, $\min_i$ and $\max_i$, respectively, of $|h-E_i|/U_i$, for
$i=1,\,2,\,3$, over all the function fields we considered of a fixed $q$ and $g$. This table provides another means to compare $E_1$ and $E_2$
and also to answer the question of which $U_i$ provides the sharper error bound. For every genus and constant field that we tested, there were
several examples for which the estimates $E_1$ and $E_2$ yielded extremely accurate estimates. In fact, there were two genus $5$ function fields
of characteristic $97$ for which $E_2 = h$.
\begin{table}[ht]
\caption{Minimum and Maximum Values of $|h-E_i|/U_i$ for $i = 1, 2, 3$}
$$\begin{array}{|rr||c|c|c||c|c|c|}
\hline
q & g  & \min_1 & \min_2 & \min_3 & \max_1 & \max_2 & \max_3 \\
\hline
997    & 3 & 0.000125 & 0.000020 & 0.000024 & 0.929573 & 0.694462 & 0.918213 \\
10009  & 3 & 0.000070 & 0.000023 & 0.000031 & 0.953275 & 0.713493 & 0.948725 \\
100003 & 3 & 0.000012 & 0.000013 & 0.000017 & 0.975737 & 0.730639 & 0.973334 \\
\hline
997    & 4 & 0.000048 & 0.000054 & 0.000068 & 0.835200 & 0.663961 & 0.825064 \\
10009  & 4 & 0.000001 & 0.000006 & 0.000008 & 0.833682 & 0.666135 & 0.831041 \\
\hline
97     & 5 & 0.000069 & 0        & 0        & 0.813107 & 0.754499 & 0.791369 \\
997    & 5 & 0.000012 & 0.000000 & 0.000000 & 0.833248 & 0.826243 & 0.839193 \\
\hline
97     & 6 & 0.000019 & 0.000072 & 0.000075 & 0.808939 & 0.789328 & 0.821493 \\
463    & 6 & 0.000032 & 0.000010 & 0.000011 & 0.755526 & 0.733616 & 0.747599 \\
\hline
19     & 7 & 0.000033 & 0.000014 & 0.000015 & 0.579609 & 0.546870 & 0.588628 \\
97     & 7 & 0.000005 & 0.000001 & 0.000001 & 0.673038 & 0.641696 & 0.664122 \\
\hline
\end{array}$$
\label{tab:minmaxalphas}
\end{table}

%For every genus and constant field that we tested, there were several examples for which the estimates $E_1$ and $E_2$ yielded extremely accurate
%estimates. In fact, there were two genus $5$ function fields of characteristic $97$ for which $E_2 = h$ and a few examples for which $|h-E_2| < 10$.
%In addition, there was a genus $5$ function field of characteristic $997$ in which the second estimate was off by $3$. In contrast, the upper bounds,
%$U_i$, were less sharp with increasing genus, but for a fixed genus, $U_i$ was generally increasingly sharp as $q$ increased. We will explain this
%behavior by considering the averages $\ah_i(q,\,g)$. However, we do note that $U_3$ was a sharper upper bound than $U_2$ consistently, so that for
%large computations, we suggest using $E_2$ as the estimate for $h$ and $U_3$ as the upper bound on $|h-E_2|$.

In Table \ref{tab:alphas}, we list average values, $\ah_i(q,\,g) = Mean_n(|h-E_i|/U_i)$, for $i=1,\,2,\,3$ ($E_2 = E_3$) and fixed $q$ and $g$,
computed from the random sampling of $n=10000$ function fields.

\begin{table}[ht]
\caption{Comparison of the $\ah_i(q,\,g)$}
$$\begin{array}{|rrr||c|c|c|}
\hline
q & g & \lambda & \ah_1(q,\,g) & \ah_2(q,\,g) & \ah_3(q,\,g)\\
\hline
997    & 3 & 1 & 0.26832306 & 0.20003340 & 0.26448274 \\
10009  & 3 & 1 & 0.27031818 & 0.20234914 & 0.26906175 \\
100003 & 3 & 1 & 0.27227076 & 0.20408453 & 0.27187490 \\
\hline
997    & 4 & 1 & 0.19223965 & 0.15306081 & 0.19019941 \\
10009  & 4 & 1 & 0.19252978 & 0.15379110 & 0.19186318 \\
\hline
97     & 5 & 2 & 0.18195632 & 0.17143328 & 0.17981087 \\
997    & 5 & 2 & 0.19188423 & 0.18894457 & 0.19190607 \\
\hline
97     & 6 & 2 & 0.15246065 & 0.14526827 & 0.15118788 \\
463    & 6 & 2 & 0.15992960 & 0.15676849 & 0.15975657 \\
\hline
19     & 7 & 2 & 0.11428348 & 0.10135344 & 0.10909269 \\
97     & 7 & 2 & 0.12684120 & 0.12176623 & 0.12602172 \\
\hline
\end{array}$$
\label{tab:alphas}
\end{table}

As with the analogous situation in hyperelliptic function fields (see Section 6 of \cite{st02b}), we assumed that the limit of the actual
averages, $\lim_{q \to \infty}\alpha_i(q,\,g) = \alpha_i(g)$, exists for each $g$. Again, we can only at best estimate what the actual limits
are, based on experimental results. Given the behavior of the $\ah_i(q,\,g)$, the data also suggest that $\alpha_i(g)$ decreases as $g$
increases, as is the case for hyperelliptic function fields \cite{st02b}. Note also that for $\lambda = 1$, the values for $\ah_2(q, \,g)$ in
Table \ref{tab:alphas} are noticeably smaller than those for $\ah_1(q, \,g)$ and $\ah_3(q, \,g)$, whereas for $\lambda = 2$, the three values
$\ah_i(a, \,g)$ match more closely for $i = 1, 2, 3$. An analogous phenomenon can be observed in Table \ref{tab:minmaxalphas}. In the next
section, we explain why this behavior is to be expected and also explain the difficulties arising in the computation of each $\alpha_i(q,\,g)$.

\subsection{Analysis of the $\alpha_i(q,\,g)$}
\label{S:alphaanalysis}

In this section, we take a closer look at the relationship between the averages $\alpha_i(q,\,g)$ and the error bounds $U_i$. In particular, we
explain the obstructions to computing the $\alpha_i(q,\,g)$ more precisely.

We have $h = \prod_{j=1}^{2g}(1-\omega_j)$, where $\omega_1, \ldots \omega_{2g} \in \C$ are the reciprocals of the zeros of the zeta function of
$K$. Write $\omega_j = \sqrt{q}e^{i\varphi_j}$, where $i$ is a fixed square root of $-1$ and each $0\leq \varphi_j < 2\pi$ for $1 \leq j \leq
2g$. It is well-known that the $\varphi_j$ can be arranged so that $\omega_j = \overline{\omega}_{j+g}$ and $\varphi_j \equiv -\varphi_{j+g}
\ppmod{2\pi}$. We may therefore order the $\varphi_j$ so that $0\leq \varphi_j \leq \pi$ for $1\leq j \leq g$. Set
\[ G_\lambda(\varphi_1, \ldots , \varphi_g) = \sum_{j=1}^{2g} e^{(\lambda+1)i\varphi_j}
    = 2 \sum_{j=1}^g \cos \big ( (\lambda+1)\varphi_j \big ) \enspace . \]
The analysis in Section 5 of \cite{ss08} shows that for large $q$, we expect that
\begin{equation} \label{eq:G}
\alpha_i(q, g) = \mbox{Mean} \left ( \frac{|h-E_i|}{U_i} \right ) \approx
    \frac{|G_\lambda(\varphi_1, \ldots , \varphi_g)| + \eta_{1, i}(q, \lambda)}{2g + \eta_{2, i}(q, \lambda)}
\end{equation}
for $i = 1, 2, 3$, where $\eta_{1, i}(q, \lambda)$ and $\eta_{2, i}(q, \lambda)$ are correction terms that depend on $i$ and vanish for $i = 1$.
Both $\eta_{1, 2}(q, \lambda) = \eta_{1, 3}(q, \lambda)$ and $\eta_{2, 3}(q, \lambda)$ also tend to zero as $q$ grows, as does $\eta_{2, 2}(q,
\lambda) = 0$ for $\lambda$ even. However, $\eta_{2, 2}(q, \lambda) = [K:\Fq(x)] - 1 = 2$ for $\lambda$ odd.

We therefore see that the averages $\alpha_i(q, g)$ essentially depend only on the distribution of the values $\varphi_j$ around the unit
circle. For example, if each $\varphi_j$ is close to either $0$ or $\pi$, then $\alpha_i(q,g) \approx 1$ (or $\approx g/(g+1)$ for $i = 2$ and
$\lambda$ odd). On the other hand, if the average of the $\varphi_j$ is close to $\pi/2$, then $\alpha_i(q,g) \approx 0$. Based on our
experimental results, it is a reasonable assumption that over all cubic function fields over a fixed base field and and of fixed genus, the
average of the $\varphi_j$ is distributed symmetrically about $\pi/2$. As the genus increases, it becomes less likely for each $\varphi_j$, for
any given function field, to be very close to either $0$ or $\pi$, thereby making it less likely for $G_\lambda(\varphi_1, \ldots , \varphi_g)$
to be large. Hence, for increasing genus, we expect a decreasingly smaller proportion of function fields with $\ah_i(q,g)$ far away from 0,
which would explain the lower values of the $\ah_i(q,\,g)$ in Table \ref{tab:alphas} with increasing genus. This also explains why the minimum
and maximum values of the $|h-E_i|/U_i$ in Table \ref{tab:minmaxalphas} generally decrease with increasing genus.

Note also that for $\lambda = 1$, the values for $\ah_2(q, \,g)$ in Table \ref{tab:alphas} are noticeably smaller than those for $\ah_1(q, \,g)$
and $\ah_3(q, \,g)$, since the denominator in the right hand side of \eqref{eq:G} is $2g+2$ for $i = 2$ and closer to $2g$ for $i = 1, 3$. For
$\lambda = 2$, this denominator is approximately $2g$ for all $i = 1,2, 3$, and thus, the three values $\ah_i(a, \,g)$ match more closely. An
analogous phenomenon can be observed in Table \ref{tab:minmaxalphas}.

If the values $\varphi_j$ ($1\leq j \leq g$) were distributed randomly in the interval $[0,\,\pi]$, over all function fields of a fixed
extension degree, genus, and base field, then precise values of the $\alpha_i(q,\,g)$ could be obtained for each $1\leq i\leq 3$. Unfortunately,
however, this is not the case, so we cannot make this assumption. In order to determine this distribution, one must know the Haar measure of a
subgroup of the symplectic group $\mbox{Sp}(2g)$ and the corresponding measure $\mu_g$. Once this measure is known, precise values of the
$\alpha_i(g)$ may be computed directly via an integral or approximated via Riemann sums. The measure $\mu_g$ has been derived for elliptic
function fields by Birch \cite{birch}, and for hyperelliptic function fields of genus $g > 1$ by Katz, Sarnak, and Weyl \cite{weyl, ks}.
Unfortunately, $\mu_g$ is not known for any function fields of degree greater than $2$. It is however conjectured that we obtain similar results
as for hyperelliptic function fields under the same assumptions. Nevertheless, determining $\alpha_i(g)$ is very difficult, so we must rely on
the approximations given in Table \ref{tab:alphas} to achieve an average running time of Algorithm \ref{alg:classnumber} that is close to
optimal. For further details, we refer the reader to Section 6 of \cite{st02b} and Section 5 of \cite{ss08}.
%Finally, we should mention that
%another way of approaching this subject is as follows. Compute the averages $\ah_i(q,\,g)$ for a large class of function fields of a given genus
%and growing values of $q$. If those averages approach a limit as $q$ increases, then use those averages heuristically in practice in order to
%improve the running time.

In the following sections, we summarize results on the application of this data to Algorithm \ref{alg:classnumber} for large examples for which no faster method is known to exist.

\subsection{Families of Curves to Consider}
\label{S:curves}

While our algorithm is a general method to compute class numbers of purely cubic function fields of characteristic at least $5$, there exist methods which work faster than ours on some special families of cubic function fields, as noted in our introduction. In this section we note which families of curves that our method works fastest on.

For polynomials $G, H \in \Fq[x]$, the curve $Y^3=GH^2$ is equivalent to $Y^3 = G^2H$. In light of this and our observation that arithmetic is faster with curves $Y^3 = GH^2$ such that $\deg(H) \leq \deg(G)$, we restricted our attention to such curves. Clearly we did not consider any cubic curves of genus 0, 1, or 2. 

Minzlaff's algorithm \cite{minzlaff} is faster than ours on superelliptic curves over a prime field and the algorithm of Castryck, Denef, and Vercauteren \cite{cdv} should run faster than ours on non-singular curves. Therefore, we did not consider any curves such that $\deg(H) = 0$. Moreover, we did not consider any curves that were equivalent to a superelliptic curve over a prime field. More specifically, if $Y^3 = D = GH^2$, $3\mid\deg(D)$ and $a$ is a root of $D\in\Fq[x]$, then there is a transformation of $Y^3=D$ to a superelliptic curve over $\Fq(a)$. To see this transformation, let $n = \deg(D)$, $t=1/(x-a)$, and $z = Y/(x-a)^{n/3}$ so that $z^3 = F(t)$, where $\deg(F) = n-2$. Multiply this through by $\sgn(F)$, the leading coefficient of $F$, and finally let $v = \sgn(F)^{n-1}z$ and $u = \sgn(F)v$, yielding the superelliptic curve $v^3 = E(u)$, where $E\in\Fq(a)[u]$ is monic and squarefree. Therefore, if $3\mid\deg(GH^2)$ (i.e. $r = 1$), then we only considered curves such that $GH$ had no linear factors. 

The following chart organizes which families of curves we considered, by genus and unit rank, based on the degrees of $G$ and $H$.

\begin{table}[ht]
\caption{Curves we Considered}
$$\begin{array}{|c||c|c||c|c|}
\hline
 & \multicolumn{2}{|c||}{r=0} & \multicolumn{2}{|c|}{r=1}\\
\hline
g & \deg(G) & \deg(H) & \deg(G) & \deg(H)\\
\hline
4 & 3 & 2 & 3 & 3  \\
\hline
\multirow{2}{*}{5}  & 4 & 2 & 5 & 2\\
 & 5 & 1 & & \\  
\hline
\multirow{2}{*}{6} & 4 & 3 & 4 & 4\\
 & 6 & 1& & \\   
\hline
\multirow{2}{*}{7} & 5 & 3 & 6 & 3\\
& 6 & 2& & \\  
\hline
\end{array}$$
\label{tab:curves}
\end{table}

\subsection{Signature $(3,1)$ Computations}
\label{S:rank0comp}

We used the estimate $E = E_2$ and the error bound $U = U_3$ and applied the values of $\ah_3(q,\,g)$, for the largest values of $q$ in Table
\ref{tab:alphas}, to the problem of computing large class numbers of purely cubic function fields of signature $(3,1)$. We will henceforth
denote this value of $\ah_3(q,\,g)$ that we use in our computations by $\ah(g)$. We computed the divisor class numbers of 
three genus $4$ and two genus $5$, $6$, and $7$ purely cubic function fields of signature $(3,1)$. We parallelized Phases 1 and 3 of each computation, using up to $64$ processors, to find class numbers up to $26$ digits. In this section, we present the results of these calculations, including timing data and the choices of certain
variables. We began with smaller examples in order to test the choices of certain parameters, in particular, the parameter $\theta$, which
regulates how often we set a kangaroo trap, and to better estimate the expected time to compute larger divisor class numbers. We list the
divisor class numbers we computed in Table \ref{tab:r0results} with corresponding statistics in Tables \ref{tab:r0g3data}, \ref{tab:r0g3data2},
and \ref{tab:r0g4data}.

The genus $4$ curves we used for the examples in this section were:
%
%$$\begin{array}{ll}
%C_1: Y^3 = x^4 +  4767220x^3 +  9719260x^2 +  9796683x +  9650320 \\
%C_2: Y^3 = x^4 + 39760243x^3 + 80354454x^2 + 40601482x + 72689039 \\
%C_3: Y^3 = x^4 + 512964174x^3 + 604076970x^2 + 208417608x + 702771176 \enspace ,
%\end{array}$$
%
%the genus $4$ curves were:
%
$$\begin{array}{ll}
%C_4: Y^3 = x^5 + 6841 x^4 + 8688 x^3 + 6670 x^2 + 5232 x + 6608 \\
%C_5: Y^3 = x^5 + 70599x^4 + 31259x^3 + 68336x^2 + 2756x + 62207 \\
%C_6: Y^3 = x^5 + 531472 x^4 + 146921 x^3 + 387330 x^2 + 602740 x + 79247\\
C_1: Y^3 = \left(x^3 + 7765x^2 + 6170x + 7834\right)\cdot \\
\quad\quad\quad \left(x^2 + 4618x + 458\right)^2 \enspace ,\\
C_2: Y^3 = \left(x^3 + 85486 x^2 + 91842 x + 21779\right)\cdot \\
\quad\quad\quad \left(x^2 + 39078 x + 54258\right)^2\enspace ,\\
C_3: Y^3 = \left(x^3 + 404647 x^2 + 836530 x + 314589\right)\cdot \\
\quad\quad\quad \left(x^2 + 945028 x + 516357\right)^2 \enspace .
\end{array}$$
The genus $5$ curves we used were:
$$\begin{array}{ll}
C_4: Y^3 = \left(x^5 + 8703 x^4 + 5098 x^3 + 1571 x^2 + 9390 x + 9945\right)x^2 \enspace ,\\
%C_4b: Y^3 = \left(x^4 + 9728 x^3 + 4683 x^2 + 1044 x + 3841 \right)\cdot
%\quad\quad\quad \left(x^2 + 9997 x + 6006 \right)^2 \enspace ,\\
C_5: Y^3 = \left(x^5 + 43583 x^4 + 40125 x^3 + 74978 x^2 + 23924 x + 38273\right)x^2 \enspace .
%C_5b: Y^3 = \left(x^4 + 81382 x^3 + 69684 x^2 + 94561 x + 70452 \right)\cdot
%\quad\quad\quad \left(x^2 + 79540 x + 94842 \right)^2 \enspace ,\\
\end{array}$$
The genus $6$ curves we used were:
$$\begin{array}{ll}
C_6: Y^3 = \left(x^4 + 212 x^3 + 980 x^2 + 939 x + 282\right)\cdot \\
\quad\quad\quad \left(x^3 + 271 x^2 + 276 x + 302\right)^2 \enspace , \\
%C_6b: Y^3 = \left(x^6 + 417 x^5 + 799 x^4 + 146 x^3 + 279 x^2 + 372 x + 951 \right)x^2 \enspace , \\
C_7: Y^3 = \left(x^4 + 4122 x^3 + 698 x^2 + 1994 x + 4252\right)\cdot \\
\quad\quad\quad \left(x^3 + 669 x^2 + 7328 x + 1019\right)^2 \enspace .
%C_7b: Y^3 = \left(x^6 + 3109 x^5 + 9366 x^4 + 591 x^3 + 8454 x^2 + 7717 x + 6227 \right)x^2 \enspace .\\
\end{array}$$
The genus $7$ curves we used were:
$$\begin{array}{ll}
C_8: Y^3 = \left(x^5 + 59 x^4 + 9 x^3 + 22 x^2 + 30 x + 37 \right)\cdot \\
\quad\quad\quad \left(x^3 + 30 x^2 + 54 x + 80 \right)^2 \enspace ,\\
%C_8b: Y^3 = \left(x^6 + 26 x^5 + 29 x^4 + 80 x^3 + 32 x^2 + 51 x + 83  \right)\cdot \\
% \quad\quad\quad \left(x^2 +89x + 83\right)^2 \enspace . \\
C_9: Y^3 = \left(x^5 + 776 x^4 + 117 x^3 + 478 x^2 + 840 x + 747 \right)\cdot \\
\quad\quad\quad \left(x^3 + 402 x^2 + 647 x + 571 \right)^2 \enspace . \\
%C_9b: Y^3 = \left(x^6 + 640 x^5 + 413 x^4 + 435 x^3 + 595 x^2 + 789 x + 427 \right)\cdot \\
% \quad\quad\quad \left(x^2 + 832 x + 22  \right)^2 \enspace . \\
%C_{10}: Y^3 = \left(x^5 + 259 x^4 + 1712 x^3 + 500 x^2 + 1156 x + 1578 \right)\cdot \\
%\quad\quad\quad \left(x^3 + 1608 x^2 + 1072 x + 1869 \right)^2 \enspace . \\
%C_{10}b:  Y^3 = \left(x^6 + 292 x^5 + 467 x^4 + 1880 x^3 + 221 x^2 + 301 x + 374 \right)\cdot \\
% \quad\quad\quad \left(x^2 + 1847 x + 443  \right)^2 \enspace . \\

\end{array}$$
In each case, we used a constant field $\Fq$ with prime $q\equiv 1\ppmod{3}$. We also have $C_i:Y^3 = G_iH_i^2$, where $G_i$ and $H_i$ are relatively prime and irreducible over the field $\Fq$ used in the respective cases; we used a random irreducible
polynomial generator supplied by NTL to choose these polynomials. In each case
The divisor class number $h$, along with the number of decimal digits in $h$ and the values $|h - E|/U$ (with $E = E_3$ and $U = U_3$), are
given for each example in Table~\ref{tab:r0results}.
%, and the number of decimal digits of each divisor class number, along with their factorizations, are given in Table \ref{tab:r0resultsb}.

\begin{table}[ht]
\caption{Divisor Class Numbers of Cubic Function Fields, Signature $(3,1)$}
$$\begin{array}{|l|l|l||c|r||c|}
\hline
\mbox{Curve} & q & g & dig. & h & |h-E|/U  \\
%\hline
%C_1    & 10^7+141 & 3 & 22 &       1000150832447729149744 & 0.0762951 \\
%C_2    & 10^8+39  & 3 & 25 &    1000018372353203578299247 & 0.2602616 \\
%C_3    & 10^9+9   & 3 & 28 & 1000020285132998304595632979 & 0.0241890 \\
\hline
%C_4 & 10^4+9     & 4 & 17 &            10226409142466713 & 0.0148396 \\
%C_5 & 10^5+3     & 4 & 20 &         99648777459613902604 & 0.0374562 \\
%C_6 & 10^6+3     & 4 & 25 &    1001264259802134080148796 & 0.3835040 \\
C_1 & 10^4 + 9 & 4 & 17 &            10011509151678732 & 0.1231612 \\ 
C_2 & 10^5 + 3 & 4 & 21 &        100380717456367838139 & 0.0009801 \\
C_3 & 10^6 + 3 & 4 & 25 &    1000964706619599167786949 & 0.1887913 \\
\hline
C_4 & 10^4 + 9 & 5 & 21 &        102398439790330982469 & 0.3211619 \\ 
%C_4b & 10^4 + 9 & 5 & 21 & 100057073626699460037 & 0.2608047 \\
C_5 & 10^5 + 3 & 5 & 26 &   10017258018358358570720475 & 0.0889572 \\
%C5b & 10^5 + 3 & 5 & 25 &  9976630448974158370762077 & 0.2349811 \\
\hline
C_6 & 10^3 + 9 & 6 & 19 &          1017494771121878691 & 0.0117266 \\
%C_6b & 10^3 + 9 & 6 & 19 & 1033829365094173632 & 0.0866109 \\
C_7 & 10^4 + 9 & 6 & 25 &    1009516362119878999248876 & 0.2450704 \\
%C_7b & 10^4 + 9 & 6 & 25 & 1004099084359934402513547 & 0.0939540 \\
\hline
C_8 & 10^2 + 3 & 7 & 15 &              117601058790012 & 0.0235252 \\
%C_8b & 10^2 + 3 & 7 & 15 & 148245067275123 & 0.0011769 \\
C_9 & 10^3 + 9 & 7 & 22 &  1002427817764983360912 & 0.1006682\\
% C_9b & 10^3 + 9 & 7 & 22 &  1052471365049681583888 & 0.1796098 \\
% C_{10} & 2011 & 7 & 24 & 126990698235970892586228 & 0.0670273 \\
% C_{10}b & 2011 & 7 & 24 & 135957648235749517073532 & 0.1451169243 \\
\hline
\end{array}$$
\label{tab:r0results}
\end{table}

%\begin{table}[ht]
%\caption{Factorization of Divisor Class Numbers, Signature $(3,1)$}
%$$\begin{array}{|l|l|l||r|l|}
%\hline
%\mbox{Curve} & q & g & dig. & \mbox{Factorization of } h \\
%\hline
%C_1 & 10^7 + 141 & 3 & 22 & 2^4 \cdot 3001 \cdot 4159 \cdot 5008303077301 \\
%C_2 & 10^8 +  39 & 3 & 25 & 19 \cdot 43 \cdot 1224012695658755909791 \\
%C_3 & 10^9 +   9 & 3 & 28 & 13 \cdot 19 \cdot 73 \cdot 114859 \cdot 482863041248304151\\
%\hline
%C_4 & 10^4 + 9 & 4 & 17 & 7 \cdot 19 \cdot 31 \cdot 1013227 \cdot 2447953\\
%C_5 & 10^5 + 3 & 4 & 20 & 2^2 \cdot  7^2 \cdot 508412129895989299\\
%C_6 & 10^6 + 3 & 4 & 25 & 2^2 \cdot  4549 \cdot 55026613530563534851\\
%\hline
%\end{array}$$
%\label{tab:r0resultsb}
%\end{table}

The divisor class number of
$\F_{10^5+3}(C_5)$ is the largest known divisor class number of a cubic function field of genus at least $4$ and signature $(3,1)$ defined by a singular curve over a large base field.

In Tables \ref{tab:r0g3data}, \ref{tab:r0g3data2}, and \ref{tab:r0g4data}, we give results from the computations of the class numbers listed in
Table \ref{tab:r0results}. Here, ``Ph.\ 1'' and ``Ph.\ 3'' give the times (in seconds) the respective phases took to complete, ``Jumps'' gives
the total number of kangaroo jumps in the computation, $\lg\theta$ indicates our choice of $\theta$, ``Traps'' records the number of
kangaroo traps that were set, $m$ is the number of processors used (if $m=1$, then a tame and a wild kangaroo ran on the same processor), and
``Total'' is the sum of these times.  ``Exp.\ 1'' gives the quantity, $\left(m|h-E|/\beta + 4\beta/m + \theta m\right)T_G$, obtained from
Proposition \ref{prop:roo0}, where $T_G$ was the time to compose two ideals in the given example, $\beta = (m/2)\sqrt{\ah(g) U}$ was the
average jump distance in the example, $E$ was the estimate of $h$, and $U$ was the upper bound on the error; the quantity Exp.\ 1 estimates the
expected time to compute the class number of the specific function field $\Fq(C_i)$ using a single processor, based on the parameters given in
Proposition~\ref{prop:roo0}. Finally, ``Exp.\ 2'' gives the quantity $\left(4\sqrt{\ah(g) U} + \theta m\right)T_G$, obtained from Proposition
\ref{prop:roo0}, which estimates the expected time to compute the divisor class number of a purely cubic function field of the given
characteristic $q$ and genus $g$ using a single processor. In Table \ref{tab:r0g4data}, we only give the total time, since Phase 1 required
very little time compared with Phase 3.

\begin{table}[ht]
\caption{Divisor Class Number Computation Data, Signature $(3,1)$, Genera $5$ and $6$}
$$\begin{array}{|l|l||r|r||r||r|r|}
\hline
\mbox{Curve} & q &\mbox{Ph.\ 1} & \mbox{Ph.\ 3} & \mbox{Jumps} & \lg\theta & \mbox{Traps}\\
%\hline
%C_1    & 10^7+141 &    553 &    4178 &    6220320 & 16 &   90 \\
%C_2    & 10^8+39  &   6098 &  114720 &  174938127 & 18 &  627 \\
%C_3    & 10^9+9   &  64577 & 1695582 & 2880612442 & 20 & 2779 \\
\hline
C_4 & 10^4 + 9 &   2727 &   73168 &   63774673 & 16 & 1087 \\
%C_4b & 10^4 + 9 & 2904 & 78832 & 63177590 & 16 & 903 \\
C_5 & 10^5 + 3 & 319184 & 3511104 & 2643608736 & 20 & 2518 \\
%C_5b & 10^5 + 3 & 343818 & 2720576 & 1925551840 & 20 & 1860 \\
\hline
C_6 & 10^3 + 9 &     25 &    5548 &    3477659 & 14 &  154 \\
%C_6b & 10^3 + 9 & 42 & 29568 & 15855545 & 16 & 300 \\
C_7 & 10^4+ 9 &   2797 & 3241472 & 1670868830 & 20 & 1516 \\
%C_7b & 10^4+ 9 & 2751 & 6147328 & 3888144966 & 20 & 3991 \\
\hline
\end{array}$$
\label{tab:r0g3data}
\end{table}
\begin{table}[ht]
\caption{Divisor Class Number Computation Data, Signature $(3,1)$, Genera $5$ and $6$}
$$\begin{array}{|l|l||r||r|r|r|}
\hline
\mbox{Curve} & q & m & \mbox{Total} & \mbox{Exp.\ 1} & \mbox{Exp.\ 2} \\
%\hline
%C_1    & 10^7+141 &  1 & 78.9\, m & 180\,  m & 243 \, m \\
%C_2    & 10^8+39  &  1 & 33.6\, h & 41.0\, h & 41.4\, h \\
%C_3    & 10^9+9   & 18 & 20.4\, d & 8.61\, d & 15.1\, d \\
\hline
C_4 & 10^4 + 9 & 16 & 21.1\, h & 14.9\, h & 11.5\, h \\
%C_4b & 10^4 + 9 & 16 & 22.7 \,h & 14.3\,h & 12.3\,h \\
C_5 & 10^5 + 3 & 64 & 44.3\, d & 25.0\, d & 32.5\, d \\
% C_5b & 10^5 + 3 & 64 & 35.5\, d & 37.8\, d & 34.5\, d\\
\hline
C_6 & 10^3 + 9 &   4 & 1.55\, h & 2.55\, h & 4.71\, h \\
%C_6b & 10^3 + 9 & 64 & 8.2\, h& 9.5\, h & 11.7\,h \\
C_7 & 10^4 + 9 & 64 & 37.5\, d & 59.0\, d & 43.4\, d \\
%C_7b & 10^4 + 9 & 64 & 71.2 \,d & 48.0\,d & 60.1\,d\\
\hline
\end{array}$$
\label{tab:r0g3data2}
\end{table}
\begin{table}[ht]
\caption{Divisor Class Number Computation Data, Signature $(3,1)$, Genera $4$ and $7$}
$$\begin{array}{|l|l|r||r|r|r||r|r|r|}
\hline
\mbox{Curve} & q & m & \mbox{Total} & \mbox{Exp.\ 1} & \mbox{Exp.\ 2} & \mbox{Jumps} & \lg\theta & \mbox{Traps} \\
\hline
%C_4 & 10^4+9     &  1 &  6.9\, m & 17.5\, m & 31.4\, m &     451594 & 14 &   27 \\
%C_5 & 10^5+3     &  2 &  9.8\, h & 10.5\, h & 16.1\, h &   35641743 & 16 &  522 \\
%C_6 & 10^6+3     & 20 & 72.4\, d & 45.2\, d & 30.3\, d & 4872971415 & 20 & 4597 \\
C_1 & 10^4 + 9 &   4 & 14.6\, m & 38.2\, m & 46.2\, m &     663789 & 14 &   45 \\ 
C_2 & 10^5 + 3 &   4 & 4.86\, h & 4.99\, h & 16.9\, h &   15961807 & 16 &  283 \\
C_3 & 10^6 + 3 &  64 & 33.6\, d & 25.9\, d & 26.5\, d & 2699578969 & 20 & 2533 \\
\hline
C_8 & 10^2+3   &  16 & 18.1\, m & 27.3\, m & 38.2\, m &    262189  & 14 &  87 \\ 
%C_8b & 10^2+3   &  16 & 78.7\, m & 30.9\,m & 52.1\, m & 2223184 & 14 & 608 \\
C_9 & 10^3 + 9 &  64 & 10.8\, d & 13.5\, d & 16.3\, d & 408090257 & 18 & 1213 \\
% C_9b & 10^3 + 9 & 64 & 17.3 \, d & 16.5\, d & 15.1\, d & 718858220 & 18 & 2255 \\
% C_{10} & 2011 & 128 & 29.0 \, d & 80.7 \, d & 104 \, d & 1061441878 & 20 & 2194 \\
% C_{10}b & 2011 & 128 & 87.7\, d & 100 \, d & 93.3 \, d & 3700228243 & 20 & 7571 \\
\hline
\end{array}$$
\label{tab:r0g4data}
\end{table}

One observation to note is the amount of variation between the actual time to compute certain divisor class numbers and the expected time, Exp.\ 1,
to compute these values using the Kangaroo method. For any given class group, the time to compute $h$ depends on the intersection of two kangaroo
paths. For a given choice of $\beta$, there are several possible choices for a set of jumps, $\{s_1,\,\ldots,\,s_{64}\}$, under the restrictions
given in Step 5 of Algorithm \ref{alg:roo0}. The number of jumps required to compute $h$ depends uniquely on this choice.
Therefore, for one set of jumps, the computation may happen to finish earlier than expected while for another set, the computation may run longer
than expected. It is impossible to know in advance how one choice of jump distances will affect the running time.

We also counted the number of useless collisions in each example. 
The computations for curves $C_1$, $C_2$, and $C_8$ yielded no collisions. 
At most, we had 4 useless collisions, for curves $C_3$, $C_4$, and $C_5$. 
A low number of useless collisions was expected, based on the results of Stein and Teske on hyperelliptic function
fields~\cite{st05}.

In the next section, we summarize the results of regulator computations in cubic function fields of signature $(1,1;1,2)$.

\subsection{Signature $(1,1; 1,2)$ Computations}
\label{S:rank1comp}

In this section, we tested the practical effectiveness of Algorithms \ref{alg:classnumber} and \ref{alg:S-regulator1} to compute the divisor
class number and extract the ideal class number and regulator of nine purely cubic function fields of signature $(1,1;1,2)$ of genera $4$ through $7$. 
We list the ideal class number $h_x$ the regulator $R_x$ and the ratio $|h-E|/U$ in Table \ref{tab:r1results}
 and data from the Kangaroo computations
in Tables \ref{tab:r1roodata} and \ref{tab:r1roodata2}. 
%The largest two examples for each genus were computed via a parallelized approach, as in
%the signature $(3,1)$ case, using up to $64$ processors. The largest divisor class numbers we computed had $28$ and $25$ decimal digits in the genus
%$3$ and $4$ cases, respectively.

The genus $4$ curves we used for the examples in this section were:
$$\begin{array}{ll}
C_{10}: Y^3 = \left( x^3 + 2833x^2 + 2425x + 5216\right)\cdot\\
\qquad\qquad\qquad\left(x^3 + 6412x^2 + 3035x + 192\right)^2 \enspace , \\
C_{11}: Y^3 = \left(x^3 + 18559x^2 + 21371x + 89569\right)\cdot\\
\qquad\qquad\qquad\left(x^3 + 1149x^2 + 83421x + 94387\right)^2 \enspace ,\\
C_{12}: Y^3 = \left(x^3 + 545795x^2 + 378803x +44676\right)\cdot\\
\qquad\qquad\qquad\left(x^3 + 736840x^2 + 529889x +983699\right)^2 \enspace .
\end{array}$$
The genus $5$ curves we used were:
$$\begin{array}{ll}
C_{13}: Y^3 = \left(x^5 + 7166 x^4 + 3769 x^3 + 7559 x^2 + 5984 x + 9826 \right)\cdot\\
\qquad\qquad\qquad \left(x^2 + 5149 x + 8000 \right)^2 \enspace ,\\
C_{14}: Y^3 = \left(x^5 + 85771 x^4 + 65270 x^3 + 5761 x^2 + 36247 x + 18059 \right)\cdot\\
\qquad\qquad\qquad \left(x^2 + 97994 x + 77903 \right)^2 \enspace .
\end{array}$$
The genus $6$ curves we used were:
$$\begin{array}{ll}
C_{15}: Y^3 = \left(x^4 + 990 x^3 + 684 x^2 + 159 x + 403\right)\cdot\\
\qquad\qquad\qquad \left(x^4 + 235 x^3 + 621 x^2 + 727 x + 49 \right)^2 \enspace , \\
C_{16}:= Y^3 = \left(x^4 + 2267 x^3 + 941 x^2 + 3751 x + 575\right)\cdot\\
\qquad\qquad\qquad \left(x^4 + 6786 x^3 + 7043 x^2 + 9857 x + 1472 \right)^2  \enspace . \\
\end{array}$$
The genus $7$ curves we used were:
$$\begin{array}{ll}
C_{17}: Y^3 = \left(x^6 + 43 x^5 + 38 x^4 + 9 x^3 + 84 x^2 + 60 x + 16\right)\cdot \\
\quad\quad\quad \left(x^3 + 53 x^2 + 106 x + 104\right)^2 \enspace , \\
C_{18}: Y^3 = \left(x^6 + 54 x^5 + 21 x^4 + 177 x^3 + 64 x^2 + 428 x + 216\right)\cdot \\
\quad\quad\quad \left(x^3 + 63 x^2 + 866 x + 687\right)^2  \enspace .\\
% C_{19}:  Y^3 = \left(x^6 + 827 x^5 + 1927 x^4 + 1872 x^3 + 1075 x^2 + 848 x + 1979 \right)\cdot \\
%\quad\quad\quad \left(x^3 + 211 x^2 + 1944 x + 767 \right)^2  \enspace .\\

\end{array}$$
In each case, we used a constant field $\Fq$, with prime $q\equiv 2\ppmod{3}$. We also have $C_i:Y^3 = G_iH_i^2$, where $G_i$ and $H_i$ are
relatively prime and irreducible over the field $\Fq$ used in the respective cases. % As with the curves generated for the signature $(3,1)$
%examples, each $G_i$ and $H_i$ was computed using a random irreducible polynomial generator supplied by NTL.

\begin{table}[ht]
\caption{Regulators and Ideal Class Numbers, Signature $(1,1; 1,2)$}
$$\begin{array}{|l|l|l||r|r||c|}
\hline
\mbox{Curve} & q & g & h_x & R_x & |h-E|/U \\
\hline
%C_{13} & 10^6 + 37 & 3 &  18 &          55561791851695519 & 0.0622306 \\
%C_{14} & 10^7 + 19 & 3 &   3 &      333335295493450981720 & 0.4660491 \\
%C_{15} & 10^8 +  7 & 3 &   3 &   333333410692036555362600 & 0.5518563 \\
%C_{16} & 10^9 +  7 & 3 &  12 & 83333335063983400511867136 & 0.0580483 \\
%\hline
C_{10} & 10^4 +  7 & 4 &  48 &            208911295254144 & 0.2552057 \\
C_{11} & 10^5 + 19 & 4 &   3 &       33359418825784135923 & 0.2460722 \\
C_{12} & 10^6 + 37 & 4 &   3 &   333383137492309549146867 & 0.1998394 \\
\hline 
C_{13} & 10^4 +  7 & 5 &   9 &       11150551526104064200 & 0.0002640 \\
C_{14} & 10^5 + 19 & 5 &   3 &  3336516656197604996052080 & 0.0005814 \\ 
\hline
C_{15} & 10^3 + 13 & 6 &   9 &         120051218369471011 & 0.0034454 \\
C_{16} & 10^4 +  7 & 6 & 108 &     9297527414155973143524 & 0.0022478 \\
\hline
C_{17} & 10^2 +  7 & 7 &  12 &             13227046636185 & 0.0015069 \\
C_{18} & 10^3 + 13 & 7 & 162 &        6765559534411953054 & 0.0082423 \\
% C_{19} & 2027 & 7 & 3 & 46901374252860881603616 & 0.0021710 \\
\hline
\end{array}$$
\label{tab:r1results}
\end{table}

%\begin{table}[ht]
%\caption{Factorization of Divisor Class Numbers, Signature $(1,1; 1,2)$}
%$$\begin{array}{|l|l|l||r|l|}
%\hline
%\mbox{Curve} & q & g & dig. & \mbox{Factorization of } h = h_xR_x \\
%\hline
%C_{13} & 10^6 + 37 & 3 & 19 & 2\cdot 3^2 \cdot 13\cdot 383\cdot 12821\cdot 870386641\\
%C_{14} & 10^7 + 19 & 3 & 22 & 2^3 \cdot 3 \cdot 5 \cdot 11 \cdot 4441 \cdot 166667 \cdot 1023524479\\
%C_{15} & 10^8 +  7 & 3 & 25 & 2^3 \cdot 3^4 \cdot 5^2 \cdot 17 \cdot 47 \cdot 293 \cdot 154321 \cdot 1708620677\\
%C_{16} & 10^9 +  7 & 3 & 28 & 2^{10} \cdot 3^2 \cdot 7 \cdot 11 \cdot 109^2 \cdot 167 \cdot 710227281795313 \\
%\hline
%C_{17} & 10^4 + 7  & 4 & 17 & 2^{11} \cdot 3^2\cdot 7\cdot 17\cdot 479\cdot 877\cdot 10883\\
%C_{18} & 10^5 + 19 & 4 & 21 & 3^2 \cdot 29 \cdot 3257 \cdot 117728460453997 \\
%C_{19} & 10^6 + 37 & 4 & 25 & 3^2 \cdot 233 \cdot 563 \cdot 847145598742455091 \\
%\hline
%C_{20} & 10^4 +  7 & 5 & 21 & 2^3 \cdot 3^2 \cdot 5^2 \cdot 7 \cdot 139 \cdot 57299853679877 \\
%C_{21} & 10^4 + 19 & 5 & 26 & 2^4 \cdot 3 \cdot 5 \cdot 313 \cdot 1499 \cdot 1667 \cdot 707887 \cdot 75328237\\
%\hline
%C_{22} & 10^3 + 13 & 6 & 19 & 3^2 \cdot 7 \cdot 241 \cdot 1097267 \cdot 64854359 \\
%C_{23} & 10^4 +  7 & 6 & 25 & 2^4 \cdot 3^5 \cdot 6594697 \cdot 39162474089897 \\
%\hline
%C_{24} & 10^2 +  7 & 7 & 15 & 2^2 \cdot 3^4 \cdot 5 \cdot 312071 \cdot 313961 \\
%C_{25} & 10^3 + 13 & 7 & 22 & 2^2 \cdot 3^5 \cdot 13^2 \cdot 61 \cdot 8191 \cdot 81727 \cdot 163393 \\
%\hline
%\end{array}$$
%\label{tab:r1resultsb}
%\end{table}

In Table \ref{tab:r1roodata}, ``BS Jumps'' and ``GS Jumps'' refer to the respective number of baby steps and giant steps computed using the
Kangaroo method in each example. In Table \ref{tab:r1roodata2}, ``Coll.'' is the number of useless collisions in the given example, ``Time''
refers to the total time taken in the computation. ``Exp.\ 1'' gives the quantity
$$\left(\frac{2m|h-E|}{\beta+ 2(\tau-1)} + \frac{2\beta}{(2\tau-1)m} + \frac{\theta m}{\tau}\right)\left(2 - \frac{1}{\tau} \right)T_G \enspace ,$$
obtained from Proposition \ref{prop:rooinf} and its proof, where $\tau$ is as given in Table \ref{tab:tau}, $\beta = m\sqrt{(2\tau-1)\ah(g) U} -
2(\tau-1)$ was the average jump distance in the example, $E$ was the estimate of $h$, $U$ was the upper bound on the error, and $T_G$ was the
time to compute a giant step in the given example; the quantity Exp.\ 1 estimates the expected time to compute the divisor class number of the
specific function field $\Fq(C_i)$ using a single processor, based on the parameters given in Proposition \ref{prop:rooinf}. ``Exp.\ 2'' gives
the quantity $\left(4\sqrt{\ah(g) U/(2\tau-1)} + \theta m/\tau\right)\left(2 - 1/\tau \right)T_G$, which estimates the expected time to compute
the divisor class number of a purely cubic function field of the given characteristic $q$ and genus $g$, using a single processor. The
remaining columns refer to the same data as in Tables \ref{tab:r0g3data} and~\ref{tab:r0g4data}. We omitted timing data on Phases 1 and 4 since
Phase 1 took under $1$ second in each case and extracting $R_x$ from $h$ in Phase $4$ took at most $6$ seconds.

\begin{table}[ht]
\caption{Regulator Computation Data, Signature $(1,1; 1,2)$}
$$\begin{array}{|l|l|r||r|r||r|r|c|}
\hline
\mbox{Curve} & q & g & \mbox{BS Jumps} & \mbox{GS Jumps} & \lg\theta & \mbox{Traps} \\
\hline
%C_{13} & 10^6 + 37 & 3 &     2690067 &    1398522 & 14 &    97 \\
%C_{14} & 10^7 + 19 & 3 &  248438065  &  129145902 & 18 &   548 \\
%C_{15} & 10^8 +  7 & 3 &  871191196  &  452808529 & 20 &   512 \\
%C_{16} & 10^9 +  7 & 3 &  2717052625 & 1412411703 & 22 &   356 \\
%\hline
C_{10} & 10^4 +  7 & 4 &     3641959 &    1168906 & 14 &    83 \\
C_{11} & 10^5 + 19 & 4 &   188417964 &   60450981 & 18 &   233 \\
C_{12} & 10^6 + 37 & 4 & 13859958890 & 4444797687 & 18 & 16888 \\
\hline
C_{13} & 10^4 +  7 & 5 &    18446263 &    4292390 & 14 &   288 \\
C_{14} & 10^5 + 19 & 5 &   736255166 &  171327866 & 18 &   711 \\
\hline
C_{15} & 10^3 + 13 & 6 &     5471177 &    1072002 & 14 &    80 \\
C_{16} & 10^4 +  7 & 6 &  1811772316 &  355124175 & 18 &  1484 \\
\hline
C_{17} & 10^2 +  7 & 7 &     2290380 &     344984 & 12 &   143 \\
C_{18} & 10^3 + 13 & 7 &   511805801 &   76010424 & 18 &   370 \\
C_{19} & 2027 & 7 & 4838512480 & 718236161 & 20 & 888 \\
\hline 
\end{array}$$
\label{tab:r1roodata}
\end{table}
\begin{table}[ht]
\caption{Regulator Computation Data, Signature $(1,1; 1,2)$}
$$\begin{array}{|l|l|r||r|r||r|r|r|}
\hline
\mbox{Curve} & q & g & m & \mbox{Coll.} & \mbox{Time} & \mbox{Exp. 1} & \mbox{Exp. 2} \\
\hline
%C_{13} & 10^6 + 37 & 3 &  1 & - &  67.3\, m &  56.6\, m &  91.8\, m \\
%C_{14} & 10^7 + 19 & 3 &  4 & 6 & 136.8\, h &  27.2\, h &  20.1\, h \\
%C_{15} & 10^8 +  7 & 3 & 42 & 2 &  20.0\, d &  10.7\, d &   7.1\, d \\
%C_{16} & 10^9 +  7 & 3 & 42 & 0 &  61.3\, d &  49.9\, d &  82.3\, d \\
%\hline
C_{10} & 10^4 +  7 & 4 &   1 & - &  41.3\, m &  60.6\, m &  52.1\, m \\
C_{11} & 10^5 + 19 & 4 &   4 & 0 &  56.9\, h &  72.1\, h &  43.4\, h \\
C_{12} & 10^6 + 37 & 4 &  64 & 6 & 383.2\, d & 125.5\, d & 128.2\, d \\
\hline
C_{13} & 10^4 +  7 & 5 &  16 & 0 &  14.9\, h &  16.3\, h &  30.9\, h \\
C_{14} & 10^5 + 19 & 5 & 64 & 0 &  32.1\, d &  48.8\, d &  89.2\, d \\
\hline
C_{15} & 10^3 + 13 & 6 &   8 & 0 &  4.83\, h &  12.9\, h &  25.2\, h \\
C_{16} & 10^4 +  7 & 6 &  64 & 1 &  65.3\, d &  89.5\, d & 176.4\, d \\
\hline
C_{17} & 10^2 +  7 & 7 &  16 & 2 &  138.\, m &  64.4\, m &  127.\, m \\ 
C_{18} & 10^3 + 13 & 7 &  64 & 0 &  25.9\, d &  27.3\, d &  51.3\, d \\
% C_{19} & 2027      & 7 & 128 & 0 & 222 \, d & 159 \, d & 313 \, d \\
\hline
\end{array}$$
\label{tab:r1roodata2}
\end{table}

\section{Conclusions and Future Work}
\label{S:conclusions} Using current implementations of the arithmetic of purely cubic function fields of signatures $(3,1)$ and $(1,1;1,2)$,
divisor class numbers up to $26$ digits were computed using the method of Scheidler and Stein \cite{ss07} and the Kangaroo algorithm as a
subroutine. In the signature $(1,1;1,2)$ case, we determined regulators up to $25$ digits. We computed approximations $\ah(q,\,g)$ of
$\alpha(q,\,g)$ for a few $q$ and for genera $3\leq g\leq 7$. This allowed us to achieve a constant-time speed-up in our computations by focusing
our effort to find $h$ on the center of the interval $[E-U,\,E+U]$, where $h$ is more likely to be found. Further speed-ups were obtained in $\R$
by computing approximations of the ratio $\tau = T_G/T_B$. The divisor class numbers are the largest such known for any cubic function field of
genus greater than $4$ constructed from a singular curve over a prime field and the regulators are the largest known regulators of any cubic function field. Moreover, the improvement to the Kangaroo
algorithm applies to the infrastructure of any signature $(1,1;1,2)$ function field.

An extension of our techniques to the case of signature $(1,1;1,1;1,1)$, with appropriate adaptations to the Baby Step-Giant Step and
Kangaroo algorithms in a two-dimensional infrastructure, is work in progress. In addition, efficient ideal and infrastructure arithmetic needs
to be developed for arbitrary (i.e.\ not necessarily purely) cubic function fields as well as for characteristic $2$ and $3$ in order to apply
this method to such function fields.
\bibliographystyle{amsplain}
\bibliography{cffcomp}
\end{document}